\documentclass{article}
\usepackage[a4paper, total={6in, 8in}]{geometry}
\usepackage{graphicx} 
\usepackage{bm}
\usepackage{amsmath}
\usepackage{amssymb}
\usepackage{mathtools}
\usepackage{url}
\usepackage{multirow}
\usepackage{adjustbox}
\usepackage{booktabs}
\usepackage{authblk}
\usepackage[english]{babel}
\usepackage{algorithm}
\usepackage{algpseudocode}
\DeclareMathOperator*{\argmax}{arg\,max}

\usepackage[bb=boondox]{mathalfa}
\usepackage{natbib}
\usepackage{subcaption}

\title{Relative Explanations for Contextual Problems with Endogenous Uncertainty: An Application to Competitive Facility Location}
\author[1]{Jasone Ramírez-Ayerbe}
\author[1]{Emma Frejinger}
\affil[1]{CIRRELT and Department of Computer Science and Operations Research, Université de Montréal, Montréal, Canada \par
\small{\texttt{miren.jasone.ramirez.ayerbe@umontreal.ca},\quad \texttt{emma.frejinger@umontreal.ca}}}
\date{}

\newcommand{\contV}{\boldsymbol{x}} 
\newcommand{\distFn}{\delta} 
\newcommand{\E}[1]{\mathbb{E}_{#1}}
\newcommand{\Prob}{\mathbb{P}}

\providecommand{\keywords}[1]
{
  \small	
  \textbf{\textit{Keywords---}} #1
}

\newtheorem{definition}{Definition}[section]
\newtheorem{example}{Example}[section]

\begin{document}

\maketitle

\section*{Abstract}

In this paper, we consider a contextual stochastic optimization problem in which unknown parameters follow distributions that depend on contextual covariates and decisions. The problem is motivated by transportation infrastructure decisions such as facility location or network design. In such high-stakes settings, decisions must often be communicated, justified, and reconsidered under alternative stakeholder requirements.  To this end, we propose a framework for computing relative counterfactual explanations. These explanations identify the smallest changes in the covariates required for a solution to satisfy prescribed constraints while limiting the performance loss to a controlled level. Whereas relative explanations have been introduced in prior literature, to the best of our knowledge this is the first work focusing on problems with binary decision variables and endogenous uncertainty.

We propose a methodology that uses the Wasserstein distance as a regularization term in the objective. Beyond improving tractability, this regularization yields explanations with desirable structural properties: it produces sparser counterfactuals, induces smoother transitions in the underlying choice distributions, and keeps the counterfactual behavior close to realistic demand patterns.  We illustrate the method using a choice-based competitive facility location problem and present numerical experiments that demonstrate its ability to efficiently compute sparse, plausible, and interpretable explanations. We further validate the framework on a real-world case study of electric vehicle charging station planning in Montreal, where the explanations reveal the minimal capacity investments and environmental conditions required to justify including a candidate location in the charging network.

\keywords{Counterfactual Explanations, Explainability in Mathematical Optimization, Contextual Optimization, Facility Location, EV charging}

\section{Introduction}

In transportation systems, demand is often shaped by the supply decisions a planner makes. Opening a charging station, locating transit stops, setting a price, or designing a service alters the alternatives users face and, in turn, the choices they make. This results in \emph{endogenous} (or decision-dependent) uncertainty, where demand distributions
depend on the decisions themselves \citep{FrejingerHewitt25}. It is commonly captured by embedding a random utility maximization (RUM) model as part of a mathematical programming formulation. There is a growing body of related literature in transportation \citep[e.g.,][]{CummingsEtAl25,BertsimasEtAl20-TS}.

In this work we focus on infrastructure planning and facility location problems in particular. In such settings, an optimal solution to a model is rarely sufficient on its own. Whether discussing the placement of healthcare facilities or an expansion of Montreal’s electric vehicle (EV) charging network, these strategic and tactical decisions involve significant capital investment and long-term societal impact \citep{chen2022using, lamontagne2025makes}. Recent qualitative work suggests that optimization models can be perceived as opaque ``black boxes'', with outputs that are not always easy to interrogate or justify \citep{lawless2025magical}. In high-stakes environments, this lack of transparency may create a barrier to adoption.

Consider the challenge of placing EV charging stations across the Island of Montreal. A planner may identify a specific station as a priority, yet it could be excluded from the solution to the optimization model. To understand this decision, we may ask: What would have to change about this specific station for it to become part of the solution? We can provide an explanation identifying the minimal, realistic modifications to a location's attributes, like adding a small number of charging outlets or improving proximity to local amenities.

Transparency is particularly challenging in contextual stochastic optimization. In these problems, the parameters, such as costs, rewards, or demand, are not fixed. Instead they follow distributions conditional on contextual covariates that describe the environment. In our Montreal example, these covariates are station- and location-specific features that shape user behavior. Examples of such covariates include the number of charging outlets available at a station or the proximity of a station to residential areas or nearby amenities \citep{lamontagne2025makes}. In addition to covariates, the probability distributions depend on the decisions themselves. In the EV context, the probability that a user chooses a station depends not only on the neighborhood's characteristics but also on which other stations the planner decides to open. This interdependency means that providing an explanation goes beyond classical sensitivity analysis.

To address these challenges formally, we consider a contextual stochastic optimization problem under endogenous uncertainty. Let $\bm{z}\in \mathcal{Z}$ denote the vector of decision variables, where $\mathcal{Z}$ is the feasible set, which may be continuous, discrete, or a mixed-integer domain; let $\bm{y}$ denote the vector of random variables capturing the endogenous uncertainty; and let $\contV^0\in \mathcal{X} \subseteq \mathbb{R}^{d_x}$ represent continuous covariates describing the context. The superscript $0$ indicates that they are the factual (i.e., original) covariates.

Uncertain parameters follow a conditional probability distribution $\Prob(\bm{y}| \bm{z}, \contV^0)$. We assume that this distribution is known, and we do not tackle the task of estimating this distribution from data \citep[see survey][]{SADANA2025271}. 

The decision-maker seeks solutions that maximize the expected reward $\E{\Prob(\bm{y}\mid \bm{z}, \contV^0)} r(\bm{y}, \bm{z})$. The resulting contextual optimization problem is given by:
\begin{equation}
\label{eq:P1}
    \bm{z}^*(\contV^0) \coloneqq \argmax_{\bm{z} \in \mathcal{Z}} \E{\Prob(\bm{y}\mid \bm{z}, \contV^0)} r(\bm{y}, \bm{z}).
\end{equation}

The literature on counterfactual explanations for optimization problems uses different definitions of explanations. We follow the definition introduced by \citet{kurtz2025counterfactual} regarding relative counterfactual explanations.  We introduce a desired feasible space  $\mathcal{D}$ representing additional stakeholder requirements not satisfied by the factual solution. We assume that $\bm{z}^*(\contV^0)\cap \mathcal{D}= \emptyset $; otherwise, the explanation problem would be trivial. 

Relative counterfactual explanations aim to identify minimal changes in covariates such that there exists a feasible solution in $\mathcal{D}$ whose expected performance is at least a prescribed fraction of the factual optimum. They are defined as follows.

\begin{definition}\label{def:relative}
\citep[Relative Counterfactual Explanations;][]{kurtz2025counterfactual} For a given factor $\alpha\in (0, \infty]$ and a desired space $\mathcal{D}$, a relative explanation for \eqref{eq:P1} is a new set of covariates $\contV$ such that there exists a feasible solution in $\mathcal{D}$ whose expected reward is at least an $\alpha$-fraction of the factual optimal expected reward.
\end{definition}

While Definition~\ref{def:relative} specifies the feasibility of a relative counterfactual explanation, additional considerations are required to characterize what constitutes a meaningful explanation in practice.

Across the counterfactual explanation literature, several desirable properties have been identified \citep{guidotti2024counterfactual}. In particular, explanations are typically required to satisfy a notion of \emph{proximity}, involving minimal changes to the factual covariates, and often exhibit \emph{sparsity}, modifying only a limited number of covariates to enhance interpretability. Other commonly discussed properties include \emph{plausibility}, which ensures that counterfactuals remain consistent with realistic data configurations, \emph{actionability}, which reflects whether the suggested changes can be implemented in practice, and \emph{robustness}, which captures the stability of explanations under small perturbations or model variations.

In contextual stochastic optimization, however, the probability distribution of the uncertain parameters is conditional on the covariates. A change in the covariates therefore induces a change in this distribution, so a feature-level dissimilarity alone may not fully capture the impact of such changes on the decision problem. This motivates the inclusion of an additional regularization term that penalizes deviations between the factual and counterfactual induced distributions, complementing the feature-based dissimilarity \citep{you2024distributional}.

In this paper, we extend the study of relative counterfactual explanations along several dimensions and make the following contributions:
\begin{itemize}
\item A methodology to compute relative counterfactual explanations for a class of contextual stochastic optimization problems under endogenous uncertainty. Existing counterfactual-explanation methods for optimization either ignore contextual covariates \citep{kurtz2025counterfactual,lefebvre2025computing,korikov2021counterfactual,korikov2023objective} or assume exogenous uncertainty \citep{forel2023explainable, viviercf}. We address the endogenous case, where the decision impacts a parametric distribution.
\item A desirable space instead of a single target solution. We generalize existing approaches that focus on explaining a single expert-provided target solution \citep{forel2023explainable,viviercf} by allowing for a desired feasible set $\mathcal{D}$ of acceptable decisions. This formulation better reflects practical settings where stakeholders impose additional constraints rather than prescribing a specific solution.
\item A tractable formulation for problems with binary decision variables, with a distributional regularizer. We combine a feature-level dissimilarity measure with a Wasserstein-based regularization on the induced distributions. This regularization yields sparser, more plausible explanations. At a moderate strength, regularization also improves computational performance relative to the unregularized formulation.

\item An instantiation on the choice-based competitive facility location problem under a multinomial logit model, as a mixed-integer bilinear program. We extend \cite{haase2009discrete} 
to a setting that captures key features of transportation infrastructure planning, where users select among facilities based on observable attributes and formulate it as a mixed-integer bilinear program.

\item A real application to EV charging network planning in Montreal. Using models estimated on revealed preference data from \citet{lamontagne2025makes}, we study two different neighborhoods and show how counterfactual explanations quantify the minimal capacity investments (additional charging outlets) and environmental conditions (proximity to nearby amenities) required to justify the inclusion of specific candidate stations in the network. To our knowledge this is the first counterfactual-explanation study for choice-based facility location, and it extends the reach of counterfactual explanations in transportation beyond the single prior shortest-path example \citep{forel2023explainable}.

\end{itemize}

The remainder of this paper is organized as follows. Section~\ref{sec:literature} discusses related work on computing counterfactual explanations for optimization problems. Section~\ref{sec:problem} introduces our problem formulation. In Section~\ref{sec:CFL}, we apply the proposed approach to a choice-based competitive facility location problem (CFLP). Computational experiments are presented in Section~\ref{sec:experiments}. We conclude in Section~\ref{sec:conclusion} and outline directions for future research.

\section{Related Work}\label{sec:literature}

We organize the related literature into three subsections. Section~\ref{sec:cf-opt} discusses counterfactual explanations for optimization problems without contextual covariate dependence. In Section~\ref{sec:cf-context}, we review counterfactual explanations for contextual optimization, where we position our methodological contribution and our transportation case study. Section~\ref{sec:cf-related} surveys related but distinct lines of work, including non-counterfactual approaches to explainability, inverse optimization, and sensitivity analysis.

\subsection{Counterfactual explanations for optimization without contextual covariates}\label{sec:cf-opt}

Within counterfactual explanations for optimization, two structural distinctions are particularly relevant for our work: the type of optimization problem (linear versus integer) and the type of  explanation considered (weak, strong, or relative).

For linear programs, \citet{kurtz2025counterfactual} formalize the three types of counterfactual explanations and propose tractable methods for the relative case by exploiting the convex structure of the underlying LP. Their illustrative examples focus on supply-chain settings such as multi-modal logistics and distribution-center reorganization. The notion of relative counterfactual that we adopt (Definition~\ref{def:relative}) follows their formalization. In a related effort, \citet{lefebvre2025computing} develop a heuristic approach for the same class of problems.

For integer linear programs, \citet{korikov2021counterfactual,korikovbeck2021counterfactual,korikov2023objective} compute weak counterfactual explanations via inverse optimization. Their setting is restricted to perturbations of cost parameters, and only those associated with variables that define the desired space. More recently, \citet{engelhardt2025counterfactual} study counterfactual explanations for integer optimization problems with mutable constraint parameters, with experiments on knapsack and the resource-constrained shortest path. By contrast, our methodology targets binary decision variables and allows changes in covariates that affect the underlying probability distribution through an endogenous random utility model.

\subsection{Counterfactual explanations for contextual optimization}\label{sec:cf-context}

When parameters in optimization models follow distributions that are conditional on contextual covariates, counterfactual explanations describe how those covariates would need to change for a different decision to emerge. \citet{forel2023explainable} initiate this line of work, assuming that parameters are estimated through weighted sample average approximations using random forests or $k$-nearest neighbors, and minimizing the changes in the covariates required for a specific expert-provided target solution to become preferable. \citet{viviercf} extend this framework to differentiable predictors such as neural networks and introduce $\varepsilon$-explanations, which relax the requirement of an explicit target by identifying the smallest change in the covariates that renders the original decision no longer sufficiently optimal. We note that the notion of relative explanation in this line evaluates both the target and the prescribed solution under the same modified covariates: a context qualifies as a relative explanation when the target decision attains a better objective than the prescribed decision in that new context. This differs from the definition we adopt (Definition~\ref{def:relative}, following \cite{kurtz2025counterfactual}), where the counterfactual solution under the modified covariates is benchmarked against the factual optimum evaluated under the original covariates.

Our work differs from this line in two respects, each motivated by our application setting. First, \citet{forel2023explainable} and \citet{viviercf} consider exogenous uncertainty, i.e., parameter distributions that are independent of the decision, whereas we explicitly address \emph{endogenous} uncertainty, which is essential in choice-based facility location since the decision to open a station reshapes the user choice probabilities. Second, rather than restricting the explanation to a single expert-provided target, we encode the set of acceptable outcomes through a desired feasible set $\mathcal{D}$, of which the expert-target case is a special instance ($\mathcal{D}$ a singleton).

Beyond methodology, the application reach of contextual counterfactual explanations remains limited. Although transportation applications such as park-and-ride and mobile-clinic deployment have motivated rich facility-location modeling \citep{HolguinVeras2012TRR, chen2022using}, the development of counterfactual explanations for such problems has not been pursued. To the best of our knowledge, the only prior work in contextual counterfactual explanation with a transportation application is the shortest-path study of \citet{forel2023explainable}; the present paper extends this reach to choice-based facility location with endogenous demand uncertainty, demonstrated on a real-world EV charging case study.

\subsection{Related but distinct lines of work}\label{sec:cf-related}

Beyond counterfactual explanations, several streams of work address the explainability of optimization solutions or build on neighboring problem formulations. We highlight here three that share themes with our setting but pursue different goals.

A first stream constructs explanations without adopting a counterfactual perspective. \citet{aigner2024framework} introduce an additional objective into the optimization problem to enhance explainability, penalizing deviations from historical solutions. While this improves interpretability, it does not provide explicit changes in the covariates needed to achieve a desired solution. \citet{goerigk2023framework} propose an interpretable model in which a decision tree partitions the cost-parameter space, with each leaf corresponding to an optimal solution for a region of that space. This approach does not generalize to contextual optimization problems with endogenous uncertainty.

A second stream is inverse optimization, and in particular contextual inverse optimization \citep{contextinverseopt}, where the covariates and the observed decisions are known but the underlying objective parameters are hidden. This is the inverse of our setting: we instead seek covariate changes that lead to a decision \emph{not} observed in the data, while satisfying a set of desired requirements.

A third stream is sensitivity analysis, which examines the range of covariate changes that leave the optimal decision unchanged and can therefore be viewed as providing a factual explanation. In contrast, our approach provides prescriptive information to stakeholders: it identifies the minimal environmental shifts required to move from a current optimal solution to a preferred alternative.

\section{Problem Statement} \label{sec:problem}
Consider the factual contextual optimization problem~\eqref{eq:P1}, and let $\bm{z}^0 \in \bm{z}^*(\contV^0)$ denote a corresponding optimal decision. Our objective is to find the smallest change in the covariates that results in a relative explanation. 

Given a function $\mathcal{L}(\contV^0,\contV)$, modeling the cost of modifying the factual covariate $\contV^0$ into a counterfactual one $\contV$, computing relative explanations for contextual optimization problems can be formulated as the following non-convex optimization problem:
\begin{subequations}
\label{eq:CE_1}
\begin{align}
\mathcal{L}^{*}_{\text{free}} \coloneqq \min_{\contV\in \mathcal{X},\, \bm{z}\in \mathcal{Z}} \quad & \mathcal{L}(\contV^0,\contV) \\
\text{s.t.} \quad & \E{\Prob(\bm{y} \mid \bm{z}, \contV)} r(\bm{y}, \bm{z}) \geq \alpha \cdot \E{\Prob(\bm{y} \mid \bm{z}^0, \contV^0)} r(\bm{y}, \bm{z}^0) \label{eq:CE_1_objConstr}\\
& \bm{z}\in \mathcal{D}  \label{eq:CE_1_DConstr}
\end{align}
\end{subequations}
where $\bm{z}$ is the decision under the modified covariates $\contV$ and $\alpha \in (0,\infty]$ is a predefined multiplicative factor that sets the required relative performance of the counterfactual solution. Note that $\bm{z}$ depends on $\contV$ (i.e., $\bm{z} = \bm{z}(\contV)$), but for readability we do not write this dependence explicitly.
We use the subscript ``free" to designate that we have not made any assumptions regarding $\Prob(\bm{y}|\bm{z},\bm{x})$.
In contrast, we denote the optimal value $\mathcal{L}^{*}_{\mathcal{H}}$ when $\Prob(\bm{y}|\bm{z},\bm{x})$ is restricted to a hypothesis class $\mathcal{H}$ of parametric distributions:
\begin{equation}
\label{eq:CE_H}
\begin{aligned}
\mathcal{L}^{*}_{\mathcal{H}} \coloneqq \min_{\contV\in \mathcal{X},\, \bm{z}\in \mathcal{Z}} \ & \mathcal{L}( \contV^0,\contV) \\
\text{s.t.}\ & \text{\eqref{eq:CE_1_objConstr}--\eqref{eq:CE_1_DConstr}, with } \Prob(\bm{y} \mid \bm{z}, \contV) \in \mathcal{H}.
\end{aligned}
\end{equation}
In Section~\ref{sec:cost}, we discuss the choice of cost function $\mathcal{L}( \contV^0,\contV)$ and introduce a Wasserstein regularization. We devote Section~\ref{sec:distributions} to a discussion of the model-free case and its use as a tractable lower bound.

\subsection{Cost Function and Wasserstein Regularization}\label{sec:cost}

A central element in computing counterfactual explanations is the cost function, which quantifies the changes between the factual covariates and their counterfactual counterparts. Classical approaches to counterfactual explanations for classification and optimization problems typically focus on minimizing dissimilarity -- using norms such as the Euclidean distance -- and complexity, i.e., the complement of sparsity, to identify a minimal set of covariates that need to be modified \citep{CARRIZOSA2024399, MILLER20191}. Sparsity can be directly enforced via the $\ell_0$ norm or induced using the $\ell_1$ norm.

We adopt a cost function $\mathcal{L}$ that combines a dissimilarity term $\mathcal{J}(\contV^0,\contV)$, which captures the magnitude of the change in the covariates, with a regularization term that we attend to later.
The dissimilarity $\mathcal{J}$ promotes proximity (small covariate changes) and, depending on the norm used, sparsity (few covariates changed); it can be measured in the factual covariate space or, when appropriate, in a transformed feature space, as the following example illustrates.

\begin{example}\label{example1}
\textbf{(Multinomial Logit Model with Availability Constraints)} Consider Problem~\eqref{eq:P1}, where the distribution $\Prob(\bm{y} \rvert \bm{z}, \contV)$ follows a multinomial logit (MNL) model, representing user preferences for charging locations. Let \( D \) denote the set of candidate charging stations, and let $z_d \in \{0,1\}$ indicate whether station $d\in D$ is opened. The covariates $\contV_d$ associated with each station $d$ may include attributes such as the number of outlets or proximity to amenities.  In this setting, the probability of a user choosing station $d$ is given by
\[
\Prob(y = d \rvert  \bm{z} ,\contV) = 
\frac{
\exp\left( {\beta}^\top \contV_d \right) z_d
}{
\sum_{k \in D} \exp\left( {\beta}^\top \contV_k \right) z_k
},
\]
where ${\beta}$ is the given vector of model parameters, and $y$ denotes the user's chosen station.

Since the MNL probabilities depend on the covariates only through the nonlinear transformation 
$\phi(\contV_d) \coloneqq \exp({\beta}^\top \contV_d)$, the model can be equivalently expressed in terms of the transformed variables $\phi(\contV_d)$. In our counterfactual formulation, we optimize directly over these transformed variables. This allows us to interpret the explanation in terms of the ``attractiveness'' of a location. Measuring dissimilarity in the $\phi$-space ensures that the explanation identifies the most significant attribute changes without the computational burden of handling exponential terms inside the optimization loop.

We define the dissimilarity as
\[
\mathcal{J}(\contV^0, \contV) = \ell_p(\phi(\contV^0), \phi(\contV)).
\]
For instance, the $\ell_1$ norm minimizes the total magnitude of the change while, as is well known, favoring vectors with few nonzero components, thereby inducing sparsity. In our context, this means the model will suggest modifying only a small number of station attributes (e.g., adding outlets to one or two specific rejected stations) rather than proposing many unrealistic changes to the entire network's attributes.

\hfill$\square$
\end{example}

In contextual stochastic optimization, decisions are evaluated through expectations taken with respect to a probability distribution that depends on the covariates. Consequently, modifying $\contV$ affects not only feature values but also the induced distribution $\Prob(\bm{y}\mid \bm{z},\contV)$ that determines the expected reward. Measuring covariates changes solely at the feature level therefore does not fully reflect their impact on the distribution underlying the objective. In settings where model outputs are distribution-dependent, recent work has advocated measuring counterfactual changes at the distributional level rather than solely in feature space \citep{you2024distributional}.

To account for this distributional effect, we introduce an additional regularization term $\Omega(\contV^0,\contV)$ that penalizes discrepancies between the distributions induced by the factual and counterfactual covariates. This leads to the following general form of the cost function:
\begin{equation}
\label{eq:costfunction}
\mathcal{L}(\contV^0, \contV) = \mathcal{J}(\contV^0, \contV) + \lambda \Omega(\contV^0, \contV),
\end{equation}
where $\lambda \geq 0$ is a hyperparameter. 
The term $\mathcal{J}$ promotes proximity and sparsity of the covariate modifications, as discussed above, while the role of $\Omega$ is to limit the deviation between the induced probability distributions.

In this work, we instantiate $\Omega(\contV^0, \contV)$ with the squared 2-Wasserstein distance $\mathcal{W}_2^2(\Prob^0, \Prob)$, motivated by its tractable finite-dimensional formulation and its established use in stochastic optimization. In discrete settings, the Wasserstein distance can be formulated as a finite-dimensional optimal transport problem \citep{PanaretosWasserstein, peyre2019computational, solomon2018optimal, villani2009optimal}:

\begin{definition}\label{def:wasserstein}
Let $\Prob^0$ and $\Prob$ be two discrete probability distributions supported on a finite set $\mathcal{S}$. The squared 2-Wasserstein distance between them is defined as
\begin{subequations}
\label{eq:wassersteindef}
\begin{align}
\mathcal{W}_2^2(\Prob^0, \Prob) \coloneqq \min_{\pi \in \Pi} \quad 
& \sum_{s \in \mathcal{S}} \sum_{s' \in \mathcal{S}} \pi_{ss'} \, \delta(s, s')^2 \\
\text{s.t.} \quad 
& \sum_{s' \in \mathcal{S}} \pi_{ss'} = \Prob^0(s) \quad \forall s \in \mathcal{S}, \\
& \sum_{s \in \mathcal{S}} \pi_{ss'} = \Prob(s') \quad \forall s' \in \mathcal{S}, \\
& \pi_{ss'} \geq 0,
\end{align}
\end{subequations}
where $\pi$ is the optimal transport plan and $\delta(s,s')$ denotes the cost of transporting a unit mass from $s$ to $s'$. When $s,s' \in \mathcal{S} \subset \mathbb{R}^n$, $\delta(s,s')$ is typically taken as the Euclidean distance.
\end{definition}

We denote by $\mathcal{L}^{*\mathcal{J}+\mathcal{W}}_{\text{free}}$ the optimal value of~\eqref{eq:CE_1} when the objective function is given by~\eqref{eq:costfunction} with the 2-Wasserstein distance, i.e.,
\begin{equation*}
\label{eq:CE_JW}
\begin{aligned}
\mathcal{L}^{*\mathcal{J}+\mathcal{W}}_{\text{free}} \coloneqq \min_{\contV \in \mathcal{X},\, \bm{z} \in \mathcal{Z}} \ & \mathcal{J}(\bm{x}^0, \contV) + \lambda\, \mathcal{W}_2^2(\Prob^0, \Prob) \\
\text{s.t.}\ & \text{\eqref{eq:CE_1_objConstr}--\eqref{eq:CE_1_DConstr}}.
\end{aligned}
\end{equation*}

Analogously, when the distribution $\Prob$ is restricted to a known family $\mathcal{H}$, we denote the corresponding optimal value by $\mathcal{L}^{*\mathcal{J}+\mathcal{W}}_{\mathcal{H}}$.

Since $\mathcal{J}$ and $\mathcal{W}_2^2$ may differ substantially in magnitude in practice, in our implementation we normalize the Wasserstein term by a reference value computed from the model-free relaxation, so that $\lambda$ admits an interpretation as a relative weight between feature-level dissimilarity and distributional regularization.

The use of the Wasserstein distance as a regularization term in stochastic optimization or in loss functions for machine learning is not new. For example, the Wasserstein distance has proven effective in ensuring stability, improving generalization under distributional shifts, and offering mechanisms for adversarial robustness and parameter tuning \citep{blanchet2019robust,gao2024wasserstein,xie2021distributionally}. We use Definition~\ref{def:wasserstein} to compute the Wasserstein distance, but other definitions may be relevant. For example, \citet{bertsimas2023optimization} propose a generalization of the Wasserstein distance in which the cost $\delta$ is replaced by the cost function of the stochastic problem itself.

\subsection{Model-free Distributions to Compute Lower Bounds}\label{sec:distributions}

The structure of the probability distribution $\Prob(\bm{y}\mid\bm{z},\contV)$ affects both the formulation and its computational tractability. We distinguish a \emph{model-based} setting from a \emph{model-free} one.

In the model-based case, we assume that the probability distribution belongs to a family of distributions $\mathcal{H}$ with parameterized transformations of the covariates $\contV$. This means that changes in the covariates influence the distribution through a structured model, and any distribution within this family must be attainable given the model. In other words, the changes must follow valid substitution patterns induced by the model.
Given that parametric distributions are often used for optimization \citep{bortolomiol2021simulation,MaiLodi20}, our main focus lies on the model-based setting.

In the model-free case, we do not impose any structural assumption on the probability distribution: rather than respecting the substitution patterns of $\mathcal{H}$, the distribution is optimized directly over the space of admissible probability measures.
As a result, the optimal model-free distribution may not correspond to any $\contV \in \mathcal{X}$ admissible under $\mathcal{H}$. However, since the model-free problem relaxes the structural constraints defining $\mathcal{H}$, it provides a valid lower bound on the objective value of the model-based formulation:
\begin{equation*}
    \mathcal{L}_{\text{free}}^* \leq \mathcal{L}_{\mathcal{H}}^*.
\end{equation*}
Intuitively, the model-free bound corresponds to a best-case scenario: it quantifies the minimum distributional shift required to satisfy the constraints, regardless of whether such a shift is attainable under $\mathcal{H}$. We use this bound to accelerate the solving of the counterfactual problem. We provide more details in the following section that is devoted to the factual problem class we aim to solve. Namely, choice-based CFLP.

\section{Application to Choice-based CFLP}\label{sec:CFL}

In choice-based CFLPs, customer behavior is not deterministically known but is captured using a RUM model (MNL used in Example~\ref{example1} is a RUM model). The probability of each customer choosing one facility or another depends on preference parameters, the optimal decision, and covariates. As such, it corresponds to a contextual optimization problem subject to endogenous uncertainty. More precisely, the problem is as follows: given a budget constraint, a firm decides which subset of candidate locations to open so as to maximize the expected captured demand. The demand capture depends on customer choices, modeled using a RUM model whose probabilities depend on the features of the locations and on the decisions, which determine the available choice set. Here we assume that the RUM model has been estimated and that the model parameters capturing customers' preferences are therefore given.

In the last two decades, a growing body of literature has proposed exact methods for solving the choice-based CFLP under the MNL model \citep{AROSVERA2013277,BENATI2002518,haase2009discrete,haasecomparison, LjubicMoreno18}. There exist different linear reformulations of the MNL that lead to mixed-integer linear problems. Specifically, \citet{haase2009discrete} exploited the proportional substitution pattern induced by the MNL model to propose an integer linear programming reformulation, which we extend to the counterfactual setting by treating the utility parameters as decision variables rather than fixed constants. This formulation can be used to solve large instances (we describe the original formulation in the appendix). 

Beyond formulations based on the MNL model, the literature has explored both richer model formulations \citep{MaiLodi20, legault2024model} of the choice-based CFLP and dedicated solution algorithms. On the modeling side, robust formulations under RUM have been developed to handle parameter uncertainty \citep{dam2022submodularity, legault2024model}, as well as variants based on alternative choice models, such as those incorporating multipurpose shopping trips \citep{mendez2023store}. 
While many applications focus on strategic or tactical planning problems, there are also operational contextual facility problems, see for example \citet{chen2022using} and \citet{HolguinVeras2012TRR}. These typically concern short-term decisions in which facility deployment or operating choices must be made in response to time-varying contextual information.

In Section~\ref{sec:explanations-cflp}, we describe the relative counterfactual explanation formulation for the choice-based CFLP under the MNL model. Section~\ref{sec:lowerandwarmstart} then introduces a tractable lower bound and a warm-start procedure to accelerate the solution process.

\subsection{Relative Counterfactual Explanations}\label{sec:explanations-cflp}

Let $D$ be a set of available locations and $E$ the set of competitive locations already occupied. We denote by $C= D \cup E$ the set of all locations. The population is considered to be composed of a finite number of customers $n\in N$, where $q_n$ denotes the demand weight of customer $n$ (e.g., the fraction of the population represented by $n$).

Following the notation in \cite{legault2024model}, the fixed attributes of the customers, such as their location, are modeled by an observed vector $\theta$, and their impact on utility is captured through model parameters. Let $\beta$ denote the vector of estimated parameters of the RUM model, which we partition as $\beta=(\beta_\theta,\beta_x)$, where $\beta_\theta$ corresponds to parameters associated with fixed attributes $\theta$, and $\beta_x$ corresponds to parameters associated with covariates. The covariates $\contV^0_d$ represent facility-level attributes that a decision-maker can modify, while $\theta$ captures fixed customer-facility attributes that cannot be altered.

We consider a factual continuous covariate denoted by $\contV^0_d$ for each facility $d\in D$. Let $\varepsilon_c$ be a random term that affects the utility of each location $c \in C$ and each customer $n\in N$, so that 
\[
u^n_c(\contV^0_c, \theta, \beta, \varepsilon_c)
=
v^n_c(\contV^0_c,\theta,\beta)
+
\varepsilon_c.
\]

In RUM models with additive error terms, utility can be decomposed as above. We further distinguish between the component of the systematic utility that depends on the facility-specific covariates (which are subject to modification in the explanation) and the component that does not depend on these covariates, so that
\[
v^n_c(\contV^0_c,\theta,\beta)
=
\bar v_c(\contV^0_c;\beta_x)
+
\hat v^n_c(\theta;\beta_\theta).
\]

For a solution $\bm {z}$, the open facilities are denoted by $C_{\bm{ z}}=D_{\bm{z}}\cup E$, where $D_{\bm{ z}}=\{d\in D \mid z_d=1\}$.

We consider an MNL model, as in Example~\ref{example1}, which implies proportional substitution patterns. 
Under the RUM framework, for each customer $n \in N$, the random variable $y_n$ denotes the facility chosen among the available set $C_{\bm z}$. The MNL model specifies the choice probabilities $\Prob_n(d \mid \bm z, \contV^0) 
= \Prob(y_n = d \mid \bm z, \contV^0).
$

In that case, the factual choice-based CFLP is formulated as 
\begin{equation}
\label{eq:simplifyCFL}
\begin{aligned}
  \max_{\bm{z}\in \mathcal{Z}} \sum_{n \in N} q_n & \sum_{d\in D_{\bm{z}}} \Prob_n(d \mid \bm{z}, \contV^0) \\
  & = \sum_{n \in N} q_n \sum_{d\in D_{\bm{z}}} \frac{\exp(v_d^n)}{\sum_{c\in C_{\bm{z}}} \exp(v_c^n)}.
\end{aligned}
\end{equation}
It is a special case of~\eqref{eq:P1}, where $\bm y=(y_1,\dots,y_{|N|})$ represents the vector of customers' discrete choices and the objective is the expected captured demand.

Let $\phi_d(\contV^0)=\exp\left [\bar v_d(\contV^0_d,\beta_x) \right]$, $\forall~d\in D$, $a_{nd}=\exp{\hat v_d^n}(\theta,\beta_\theta),$ $\forall~d\in D$, $n\in N$,  and $b_n= \sum_{e\in E} b_{ne}= \sum_{e\in E} \exp{v_e^n}$, $\forall~ n\in N$. Note that $\phi_d(\contV^0)$ captures only the component of the utility that depends on the covariates of facility $d$, which are the variables we allow to change in the counterfactual explanation. 
Examples of such covariates include facility-specific attributes such as additional amenities, service quality indicators, or marketing effort, which can be modified by the decision-maker. Customer-specific effects, such as distances between customers and facilities or other fixed demographic characteristics, are incorporated in $\hat v_c^n(\theta,\beta_\theta)$ and therefore in $a_{nd}$. 
Hence, $\phi_d(\contV^0)$ does not vary across customers, whereas $a_{nd}$ does.

The probabilities that a customer $n$ chooses location $d\in D$, or $e\in E$  are respectively given by
\begin{subequations}
\label{eq:probabilities}
\begin{align}
    &\Prob_n(d| \bm{z}, \contV^0)= \frac{\phi_d(\contV^0) a_{nd} z_d}{\sum_{d'\in D} \phi_{d'}(\contV^0) a_{nd'}z_{d'}+b_n}, \text{and}\\
    & \Prob_n(e| \bm{z}, \contV^0)= \frac{b_{ne}}{\sum_{d'\in D} \phi_{d'}(\contV^0) a_{nd'}z_{d'}+b_n}.
\end{align}
\end{subequations}
Together, \eqref{eq:probabilities} specify the choice distribution over all alternatives $c\in C =D \cup E$. Let $\Prob^0_n=\Prob^0_n(c| \bm{z}^0, \contV^0)$ denote this distribution evaluated at the factual solution, for each customer $ n\in N$.

We consider the relative counterfactual explanation with Wasserstein regularization, where the cost function \(\mathcal{J}(\contV^0, \contV)\) is given by the \(\ell_1\)-norm of the transformed covariates, as in Example \ref{example1}:
\[
\mathcal{J}(\contV^0, \contV) = \sum_{d \in D} \left| \phi_d(\contV) - \phi_d(\contV^0) \right|.
\]

In the CFLP setting, we take the support of the choice distributions to be the full set $C = D \cup E$, including stations with zero choice probability. Using the fixed set $C$ ensures that the factual and counterfactual distributions share a common support, so the optimal transport plan $\pi^n_{cc'}$ is well defined. We define the transport cost $\delta(c,c')$ as the Euclidean distance 
between the geographic coordinates of facilities $c$ and $c'$. By using Euclidean distance as the transport cost $\delta(c,c')$, the Wasserstein regularization favors counterfactuals where demand is captured by nearby geographic alternatives. This prevents the model from suggesting unrealistic behavioral shifts where there is substitution with stations that are far (in distance) from the choice in the factual setting.

The optimal value $\mathcal{L}_{\mathcal{H}}^{*\mathcal{J}+\mathcal{W}}$ is obtained by solving the following non-convex optimization problem:
\begin{subequations}
\label{eq:RE1}
\begin{align}
\min_{\pi^n, \phi(\contV), \bm{z}}
& \sum_{d\in D} |\phi_d(\contV)-\phi_d(\contV^0)| + \lambda \sum_{n\in N}\sum_{c,c'\in C} \pi^n_{cc'} \distFn(c,c')^2 \label{eq:c1} \\[6pt]
\text{s.t.}\quad
& \sum_{c'\in C} \pi^n_{cc'} = \Prob^0_n(c), \qquad c\in C, \ n\in N \label{eq:w1} \\[6pt]
& \sum_{c\in C} \pi^n_{cd} = \frac{\phi_d(\contV) a_{nd} z_{d}}{\sum_{k\in D} \phi_k(\contV) a_{nk}z_{k}+b_n}, \qquad d\in D, \ n \in N \label{eq:w2} \\[6pt]
& \sum_{c\in C} \pi^n_{ce} = \frac{b_{ne}}{\sum_{k\in D} \phi_k(\contV) a_{nk}z_{k}+b_n}, \qquad e\in E, \ n \in N \label{eq:w3} \\[6pt]
& \sum_{n\in N} q_n \sum_{d\in D} \frac{\phi_d(\contV) a_{nd} z_d}{\sum_{k\in D} \phi_k(\contV) a_{nk}z_{k}+b_n} \ge \alpha \sum_{n\in N}q_n \sum_{d\in D} \Prob_n^{0} (d) \label{eq:relative_cap} \\[6pt]
& \sum_{d\in D} z_d = B \label{eq:budget1} \\[6pt]
& \bm{z} \in \mathcal{D} \label{eq:dspace} \\[6pt]
& \pi^n_{cc'} \ge 0, \qquad c,c' \in C, \ n \in N \label{eq:domain1} \\[6pt]
& z_d \in \{0,1\}, \qquad d \in D \label{eq:domain2} \\[6pt]
& \phi_d(\contV) \in \Phi, \qquad d \in D \label{eq:domainf}
\end{align}
\end{subequations}
Objective~\eqref{eq:c1} follows the general structure introduced in \eqref{eq:costfunction}. The first term corresponds to the $\ell_1$ dissimilarity measure and promotes sparse deviations of the transformed covariates $\phi(\contV)$ from their factual values. The second term represents the Wasserstein regularization, expressed through the optimal transport variables $\pi^n$, and penalizes discrepancies between the factual and counterfactual choice distributions.
Constraints~\eqref{eq:w1}-\eqref{eq:w3} impose the Wasserstein distance by ensuring a valid transport plan. Constraints~\eqref{eq:relative_cap} guarantee that the captured demand is at least an $\alpha$-fraction of the factual captured demand. Constraint~\eqref{eq:budget1} defines the budget, and Constraint~\eqref{eq:dspace} ensures that the solution $\bm{z}$ lies in the desired space. We assume that the desired space corresponds to configurations in which a subset of the binary variables $z_d$ are fixed to 1. Constraint~\eqref{eq:domainf} restricts the transformed covariates to their feasible set, i.e., the values of $\phi(\bm x)$ attainable by a feasible covariate $\bm x \in \mathcal{X}$. This set encodes, in the transformed space, any restrictions imposed on the covariates (for instance, bounds on $\bm x$).

Feasibility of \eqref{eq:RE1} depends on the attainable range of the transformed covariates $\phi_d(\contV)$. In particular, if $\phi_d(\contV)$ is unbounded from above (e.g., when $\mathcal{X}=\mathbb{R}^{d_x}_+$ and $\bar v_d$ is linear), then for any feasible $\bm{z}$ satisfying the budget constraint, the relative performance constraint \eqref{eq:relative_cap} can always be satisfied by increasing $\phi_d(\contV)$ sufficiently. Moreover, for a given $\bm{z}$ and $\contV$, the Wasserstein distance can always be computed. If the covariate domain $\mathcal{X}$ is bounded, feasibility depends on whether the target factor $\alpha$ lies within the achievable performance range.

We follow closely the reformulation in \citet{haase2009discrete} to obtain a mixed integer bilinear reformulation. Compared to the original reformulation, we add constraints to ensure the correct definition of the MNL model, given that in our case the covariates are variables rather than fixed parameters. Let $w_{nd}$ and $\ u_{ne}$ be non-negative decision variables, then calculating relative counterfactual explanations is equivalent to solving the following bilinear optimization problem: 

\begin{subequations}
\label{eq:RE2}
\begin{align}
\min_{\substack{\pi^n, \phi(\contV), \bm{z} \\ w, u}} \quad 
& \sum_{d\in D} |\phi_d(\contV)-\phi_d(\contV^0)| +\lambda \sum_{n\in N}\sum_{c,c'\in C} \pi^n_{cc'} \distFn(c,c')^2 \label{eq:REW} \\[6pt]
\text{s.t.} \quad 
& \sum_{c'\in C} \pi^n_{cc'} = \Prob^0_n(c), \qquad c\in C, \ n\in N \label{eq:W1} \\[6pt]
& \sum_{c\in C} \pi^n_{cd}=w_{nd}, \qquad d\in D, \ n\in N \label{eq:W2} \\[6pt]
& \sum_{c\in C} \pi^n_{ce}=u_{ne}, \qquad e\in E, \ n\in N \label{eq:W3} \\[6pt]
& \sum_{n\in N}q_n \sum_{d\in D}w_{nd} \ge \alpha \sum_{n\in N}q_n \sum_{d\in D} \Prob_n^{0} (d) \label{eq:relative_cap2} \\[6pt]
& \sum_{e\in E} {u_{ne}}+ \sum_{d\in D}{w_{nd}} \le 1, \qquad n \in N \label{eq:prob_valid} \\[6pt]
& \phi_{d}(\contV) a_{nd}(w_{nd}-z_d) + b_n w_{nd} \le 0, \qquad d\in D, \ n \in N \label{eq:defMNL1} \\[6pt]
& w_{nd} - \frac{\phi_d(\contV) a_{nd}}{b_{ne}} u_{ne} \le 0, \qquad d \in D, \ e \in E, \ n \in N \label{eq:defMNL2} \\[6pt]
& \frac{\phi_d(\contV) a_{nd}}{b_{ne}}u_{ne} \le w_{nd}+(1-z_d), \qquad d\in D, \ e\in E, \ n \in N \label{eq:defMNL3} \\[6pt]
 & w_{nd} \geq 0, \qquad d \in D, \ n \in N \label{eq:domainw}\\
& u_{ne} \geq 0, \qquad e \in E, \ n \in N \label{eq:domainu}\\
&\eqref{eq:budget1}\text{--}\eqref{eq:domainf}. \nonumber
\end{align}
\end{subequations}
Objective~\eqref{eq:REW} coincides with that of \eqref{eq:c1} and preserves the same interpretation: the $\ell_1$ term promotes sparse deviations of the transformed covariates, while the Wasserstein term penalizes discrepancies between the factual and counterfactual choice distributions. Constraints~\eqref{eq:W1}-\eqref{eq:W3} define the Wasserstein distance, while constraint~\eqref{eq:relative_cap2} guarantees that the captured demand is at least an $\alpha$-fraction of the factual captured demand. Constraints~\eqref{eq:prob_valid} ensure valid distributions for each customer. Constraints~\eqref{eq:defMNL1}-\eqref{eq:defMNL3} define the MNL model. Constraints~\eqref{eq:defMNL3} are the aforementioned additional constraints compared to \citet{haase2009discrete}. Unlike traditional CFLPs where the utilities are constants, our counterfactual framework treats the attractiveness $\phi_d(\contV)$ as a decision variable. These constraints are therefore necessary to respect the proportional substitution patterns implied by the MNL model.

When solving~\eqref{eq:RE2}, we directly obtain the transformed covariate quantities
$
\phi_d(\contV) = \exp\left [ \bar v_d(\contV_d;\beta_x)\right ]
$. In the case of a single covariate per facility, i.e., $\contV = (x_1, \dots, x_D)$ with each $x_d \in \mathbb{R}_+$ and $\bar v_d(x_d;\beta_x) = \beta_{d} x_d$,
we recover the feature as $x_d = \frac{1}{\beta_{d}} \log\!\big(\phi_d(\contV)\big)$.

In the case of multiple covariates, $\contV_d \in \mathbb{R}^J$, we recover them by solving, for each facility $d$,
\[
\min_{\contV_d} \;\|\contV_d - \contV_d^{\,0}\|_p
\quad\text{s.t.}\quad
\bar v_d(\contV_d;\beta_x)
=
\log\!\big(\phi_d(\contV)\big),
\]
where $\contV^0$ is the factual covariate.
Here $\|\cdot\|_p$ denotes a general $\ell_p$ norm (e.g., $\ell^2_2$ or $\ell_1$). In many applications, the contextual utility component is linear,
$
\bar v_d(\contV_d;\beta_x)
=
\sum_{j=1}^J \beta_{dj} x_{dj},
$
where $\beta_{dj}$ are elements of $\beta_x$. In this case, the post-processing problem is convex and reduces to a quadratic program with linear constraints when $p=2$, or to a linear program when $p=1$.

\subsection{Lower Bound and Warm Start} \label{sec:lowerandwarmstart}

As stated in Section~\ref{sec:distributions}, to speed up the overall solving process, we compute a lower bound $\mathcal{L}^{*\mathcal{W}}_{\text{free}}$ using the model-free formulation. 
More precisely, we solve the following mixed integer linear program: 
\begin{subequations}
\label{eq:LB}
\begin{align}
\mathcal{L}^{*\mathcal{W}}_{\text{free}} = \min_{\pi^n, w, u, \bm{z}, \bar{\bm{z}}} \quad & \sum_n\sum_{c,c'\in C} \pi^n_{cc'} \distFn(c,c')^2 \label{eq:LB_obj} \\[6pt]
\text{s.t.} \quad & w_{nd} \le \bar{z}_{nd}, \qquad d \in D, \ n\in N \label{eq:aux1} \\[6pt]
& w_{nd} \ge \epsilon \bar{z}_{nd}, \qquad d\in D, \ n\in N \label{eq:aux2} \\[6pt]
& \bar{z}_{nd} \le z_d, \qquad d \in D, \ n\in N \label{eq:aux_z1} \\[6pt]
& z_d \le \sum_{n\in N} \bar{z}_{nd}, \qquad d\in D \label{eq:aux3} \\[6pt]
& \bar{z}_{nd} \in \{0,1\}, \qquad d\in D, \ n\in N \label{eq:aux_z2} \\[6pt]
& \eqref{eq:budget1}\text{--}\eqref{eq:domain2}, \ \eqref{eq:W1}\text{--}\eqref{eq:prob_valid}, \ \eqref{eq:domainw}\text{--}\eqref{eq:domainu}. \label{eq:modelfree_f}
\end{align}
\end{subequations}
We introduce binary variables $\bar z_{nd},~\forall d\in D, n\in N$, which take the value 1 if customer $n$ has a non-zero probability of going to facility $d$, and 0 otherwise. The parameter $\epsilon > 0$ enforces that if customer $n$ is assigned to facility $d$ (i.e., $\bar z_{nd} = 1$), then the probability $w_{nd}$ must be strictly greater than zero. Constraints~\eqref{eq:aux1} and~\eqref{eq:aux2} ensure the correct definition of $\bar z_{nd}$, while Constraints~\eqref{eq:aux3} guarantee that if a facility $d$ is not opened, all customers have a zero probability of choosing it; conversely, if $d$ is opened, at least one customer must have a strictly positive probability of doing so.

Problem~\eqref{eq:LB} is a relaxation of~\eqref{eq:RE2}. It is obtained by removing the transformed covariates $\phi(\contV)$ and the MNL-defining constraints~\eqref{eq:defMNL1}--\eqref{eq:defMNL3}, while keeping the constraints related to transport plan~\eqref{eq:W1}--\eqref{eq:W3}, relative-capture~\eqref{eq:relative_cap2}, valid distribution~\eqref{eq:prob_valid}, budget~\eqref{eq:budget1}, and desirable space~\eqref{eq:dspace}. Removing the structural MNL constraints enlarges the feasible region in the probability variables $(w,u)$, so any feasible solution of~\eqref{eq:RE2} induces a feasible solution of~\eqref{eq:LB} with the same transport cost, provided $\epsilon$ in \eqref{eq:aux2} is no larger than the smallest positive choice probability attainable under~\eqref{eq:RE2}. Since~\eqref{eq:LB} minimizes only the Wasserstein term, whereas~\eqref{eq:RE2} minimizes $\mathcal{J}+\lambda\,\mathcal{W}_2^2$ with $\mathcal{J}\ge 0$, the optimal value $\mathcal{L}^{*\mathcal{W}}_{\text{free}}$ does not bound the full objective of~\eqref{eq:RE2} directly. Instead, using $\mathcal{J}\ge 0$,
\begin{equation*}
\mathcal{L}^{*\mathcal{J}+\mathcal{W}}_{\mathcal{H}} = \mathcal{J}^* + \lambda\,\mathcal{W}_2^{2*} \;\ge\; \lambda\,\mathcal{L}^{*\mathcal{W}}_{\text{free}},
\end{equation*}
so that $\lambda\,\mathcal{L}^{*\mathcal{W}}_{\text{free}}$ is a valid lower bound on the optimal value of~\eqref{eq:RE2}.

In addition, we propose the following constructive greedy algorithm to generate a feasible initial solution that can be used as a warm start:
\begin{enumerate} 
\item Fix $z_d=1$ for $d$ such that $\bm{z}\in \mathcal{D}$, otherwise set $z_d=0$. 
\item Compute the factual captured demand $ Q^{\text{factual}}=\sum_{n}q_n \sum_d \Prob_n^{0}(d). $ 
\item Sort facilities according to decreasing factual captured demand and fix $z_d=1$ for the first $B-1$ facilities in the list. 
\item Set $\bar \phi_{d}(\contV)=c$ for all $d$ such that $\bm{z}\in \mathcal{D}$ and $\bar \phi_{d}(\contV)=1$ otherwise. 
\item Increase $c$ until the relative performance constraint is satisfied, i.e., 
\[ \sum_{n \in N} \frac{\sum_{d \in D}\bar \phi_{d}(\contV) a_{nd} z_d} {\sum_{d \in D} \bar\phi_{d}(\contV) a_{nd} z_d + b_n} \geq \alpha Q^{\text{factual}}, \] 
whenever such a value exists within the admissible range of $\phi_d(\contV)$. 
\item Compute probabilities $\Prob_n(d)$ using~\eqref{eq:probabilities} with $\bar \phi_{d}(\contV)$ and evaluate the Wasserstein distance $\mathcal{W}_2^2(\Prob^0,\Prob)$. \end{enumerate}

This procedure allows us to obtain a feasible solution that respects both the demand model and the desired space constraints, while remaining computationally inexpensive.

\section{Computational Experiments}\label{sec:experiments}

This section is organized as follows. In Section~\ref{sec:illustrative}, we present a small illustrative example. In Section~\ref{sec:larger}, we assess the scalability of our approach across different problem dimensions and show that relative explanations can be computed more efficiently when the Wasserstein regularization is included. Finally, in Section~\ref{sec:casestudy}, we apply the framework to a real-world case study on EV charging station planning in Montreal, using revealed preference data from \citet{lamontagne2025makes}.

All numerical experiments were conducted on an Intel Core i9-10980XE processor. We used Python 3.11.7, and all optimization models were solved using Gurobi 11.0.1~\citep{gurobi}. The code is available at \url{https://github.com/jasoneramirez/relative-explanations-cflp}.

\subsection{Illustrative Example}\label{sec:illustrative}

We use a small example to serve two purposes. First, to illustrate a relative counterfactual explanation in the case of a choice-based CFLP. That is, the smallest change in the covariates such that a feasible configuration including a currently non-selected facility maintains the captured demand of the factual solution (i.e., $\alpha=1$). Second, to illustrate the effect of the Wasserstein regularization, which trades a larger $\ell_1$ cost for explanations whose induced choice distributions stay closer to the factual ones.

\begin{figure}[t]
\centering
\includegraphics[width=0.5\linewidth]{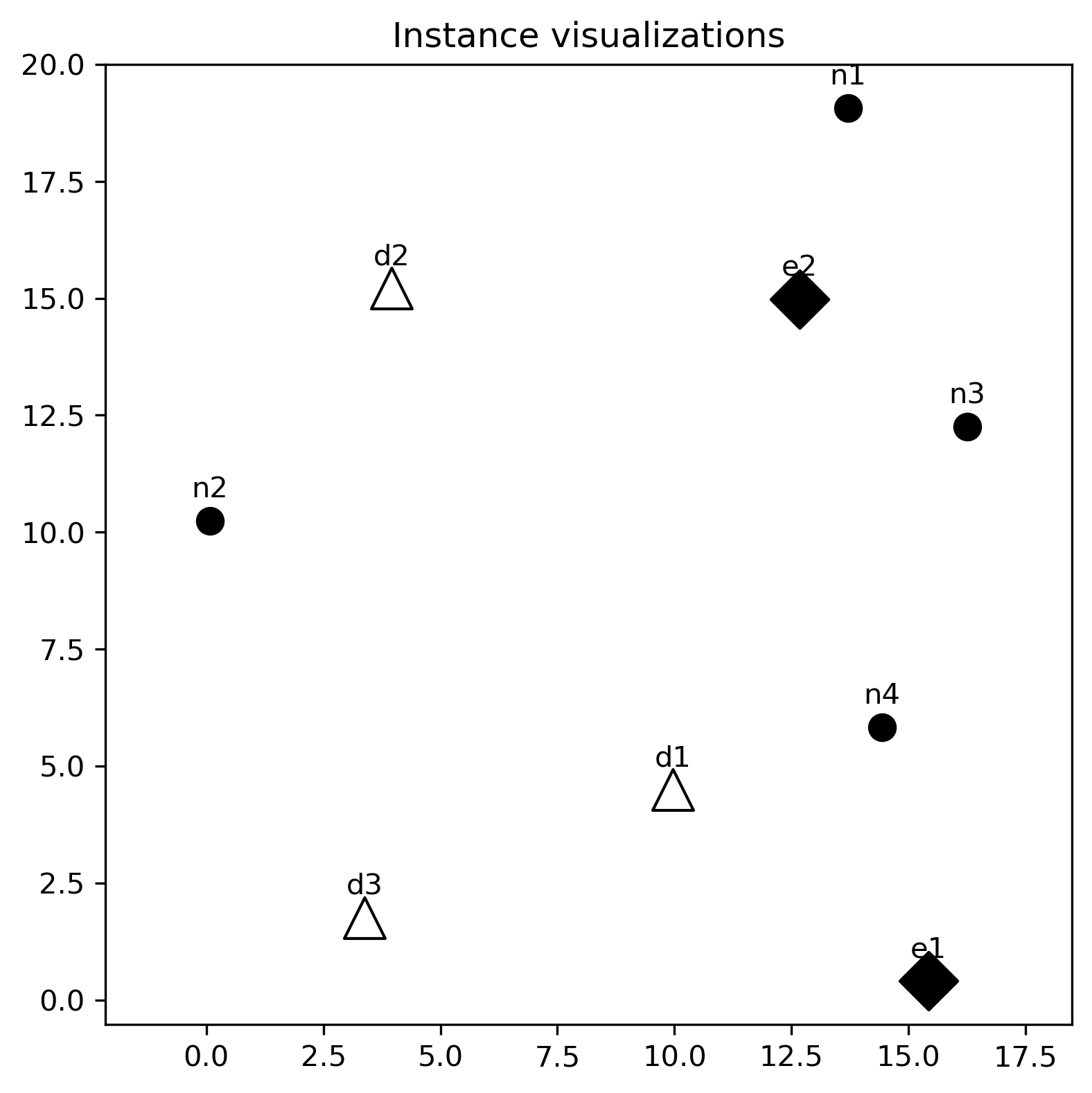}
\caption{Candidate (triangles), competitor (solid diamonds) and customer (solid circles) locations in the illustrative example}
\label{fig:illustration}
\end{figure}

We consider $|N|=4$ users, $|D|=3$ candidate facilities, $|E|=2$ competitive facilities, and $B=2$ facilities to be opened. The coordinates of $n\in N$ and $c\in C=D\cup E$ are randomly generated from a uniform distribution over the interval $[0,20]$, using random seed 10. All customers have unit demand, i.e., $q_n=1$ for all $n\in N$. Let $\theta_c^n$ be the Manhattan ($\ell_1$) distance between customer $n$ and facility $c\in D\cup E$. Let $\contV_d$ be an intrinsic attractiveness score associated with facility $d$, e.g., perceived quality, service variety, or customer satisfaction. The utility function takes the following form $
v_d^n = \hat v_d^n + \bar v_d = -0.1\,\theta_d^n + \contV_d $
for each customer $n$ and facility $d$,
and $ v_e^n = -0.1\,\theta_e^n $
for each customer and competitive facility $e$.
We assume that the factual $\contV^0_d = 0$, $\forall d\in D$. Choice probabilities follow a MNL model. The realization used in this section is depicted in Figure~\ref{fig:illustration}.

When solving \eqref{eq:simplifyCFL}, the optimal (factual) solution is to open facility 1 and 2, leading to a captured demand of $1.8804$.

Our desired space is to open facility 3. Accordingly, our relative counterfactual explanation provides the minimal changes required in $\contV_d$ such that the captured demand remains the same (with $\alpha=1$) with facility 3 opened. To obtain such an explanation, we solve~\eqref{eq:RE2}. First, we consider no regularization term, i.e., $\lambda=0$, minimizing only the changes in the covariates.

We obtain the following solution: 
\[
\contV_{d_1} = 0.350 \hspace{1cm} \contV_{d_2} = 0.000 \hspace{1cm} \contV_{d_3} = 0.000
\]
with a Wasserstein distance of 164.917 and a captured demand of 1.8804. Since $\phi_d = \exp(\contV_d)$ and the factual values satisfy $\phi_d = 1$, this corresponds to an $\ell_1$ cost of $0.419$. In other words, the attractiveness of facility~1 must increase by approximately $42\%$ in exponential scale to compensate for opening facility~3 instead of facility~2. The model increases the attractiveness of location~1, which makes it optimal to open location~3 instead of location~2.

Next, we set $\lambda=0.25$, thereby adding the Wasserstein regularization. This introduces a trade-off that encourages the new choice distributions to remain close to the factual ones while still minimizing covariate changes. The solution in this case is:
\[
\contV_{d_1} = 0.000 \hspace{1cm} \contV_{d_2} = 0.479 \hspace{1cm} \contV_{d_3} = 0.000
\]
with a Wasserstein distance of 90.8493 and a captured demand of 1.8804. Here, the corresponding $\ell_1$ cost equals $0.614$, which is larger than in the unregularized case. Hence, keeping the induced choice distributions closer to the factual ones requires a stronger increase in covariate attractiveness. Unlike the previous case, the optimal explanation now consists of increasing the attractiveness of location~2.
\begin{figure}[t]
\centering
\includegraphics[width=0.8\linewidth]{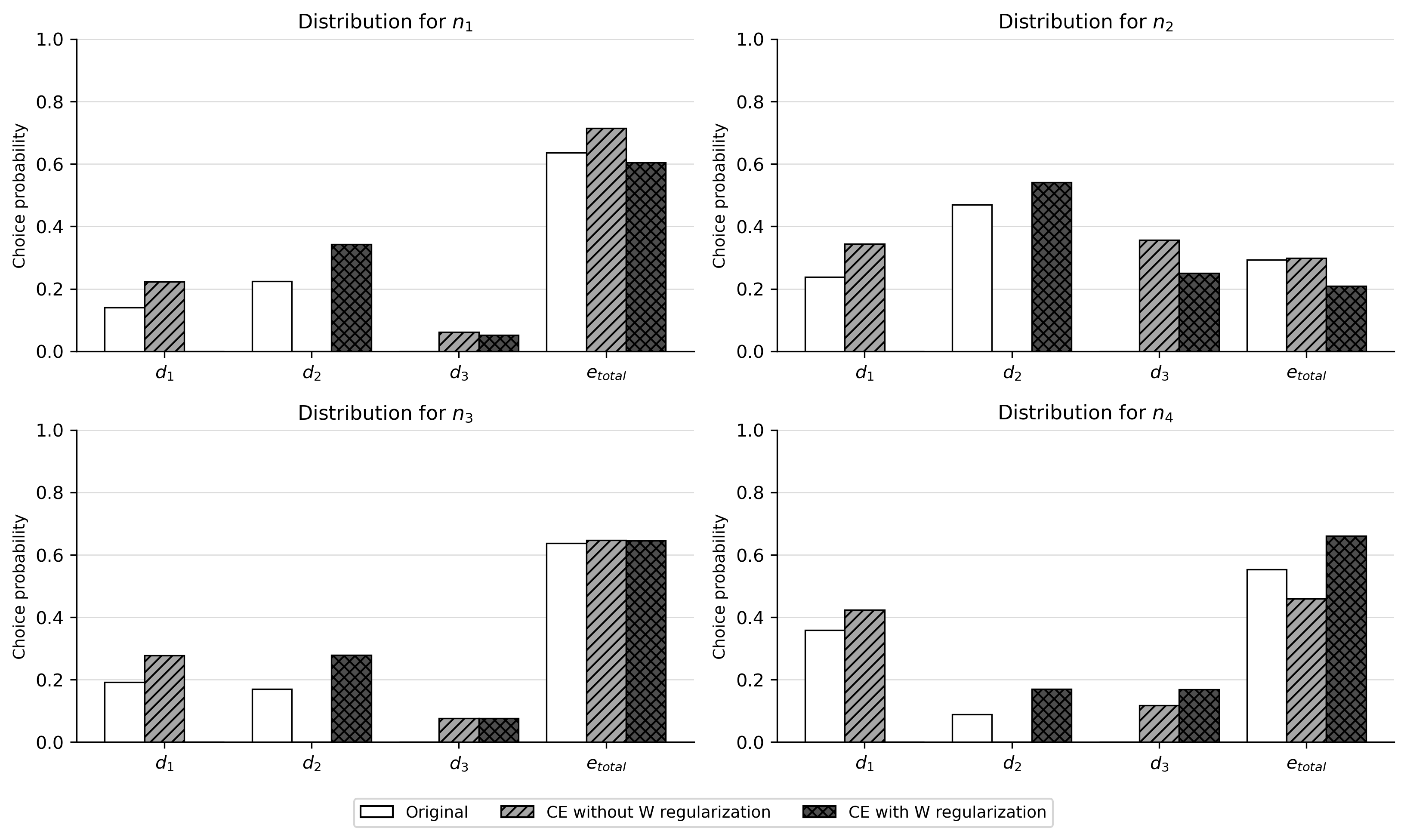}
\caption{Probability distributions for each customer in different solutions. One plot per customer. Bars show the choice probability distribution induced by each of the three solutions: the factual, and the counterfactual explanations (CE) without and with Wasserstein regularization. Competitor locations are aggregated into $e_{total}$.}
\label{fig:distributions}
\end{figure}

We show the choice probability distributions induced by the different solutions in Figure~\ref{fig:distributions}. There is one plot for each customer. Notice how, even though the covariate changes are larger in the solution obtained with regularization, the sparsity of the solution remains the same. As expected, the probability distributions are closer to the factual ones than those induced by the solution without regularization. When increasing the attractiveness of facility 2, we observe that the probability mass of the customers shifts from facility 1 to facility 3. These changes also affect the probability of choosing the competitor location.

\subsection{Larger Instances}\label{sec:larger}
In this section, we build on the illustrative example and report results for instances of different sizes and characteristics. More precisely, we vary the number of customers, the number and location of facilities, and the budget.
All facilities (candidate and competitive) and customers are independently and uniformly generated over a square region 
$[0,20]\times[0,20]$. Distances used in the utility function correspond to the Manhattan metric between customer and facility locations. We use the same utility function as before, given by $v_d^n=\hat v_d^n+\bar v_d=-0.1 \theta_d^n+\contV_d$ for each customer $n$ and facility $d$ and $v_e^n=-0.1\theta_e^n$ for each customer and competitive facility $e$. We assume that the factual $\contV^0_d=0$, $\forall d\in D$ and that all customers have unit demand, i.e.,  $q_n=1$ for all $n\in N$. Choice probabilities follow a MNL model. We solve~\eqref{eq:RE2} with $\alpha=1$ for regularization parameters $\lambda\in\{0,0.1,1\}$. Since we do not impose upper bounds on $\contV$ in these synthetic experiments, by the feasibility argument in Section~\ref{sec:explanations-cflp} a relative explanation always exists. We also measure the characteristics of the explanations, including the factual captured demand $Q^{\text{factual}}$, the new captured demand $Q^{\text{new}}$, the Wasserstein distance between the counterfactual and factual instances $\mathcal{W}_2^2$, and the complexity, defined as the percentage of covariates changed (the complement of sparsity). For each set of characteristics, we randomly generate ten instances. The computational time limit is set to one hour. For each instance, we use a solution obtained with the heuristic described in Section~\ref{sec:lowerandwarmstart} as a warm start. The time required to generate this solution is excluded from the reported computational time, which refers only to solving the counterfactual problem. The construction of the initial feasible solution takes approximately 0.1–0.5 seconds on average, while computing the Wasserstein bound requires between 0.5 and 4 seconds, depending on the problem size (see Table~\ref{tab:initial_times}). We set $\epsilon=10^{-5}$ in all our experiments.
\begin{table}[t]
\centering
\caption{Experimental results for different problem scales and regularization strengths}
\label{tab:results}
\begin{adjustbox}{max width=\textwidth}
\setlength{\tabcolsep}{4pt}
\begin{tabular}{ccccc rrrrrrrrr}
\toprule
N & D & B & $\lambda$ & $Q^{\text{factual}}$ & $Q^{\text{new}}$ & $\mathcal{W}_2^2$ & $\ell_1$ & Complexity & Avg Time [s] & Median Time [s] & TL & Gap \\
\midrule
\multirow[t]{12}{*}{100}  
  & \multirow[t]{6}{*}{10}
    & \multirow[t]{3}{*}{4} 
       & 0    & 47.612 & 47.612 & 649.462 & 0.157 & 0.330  & 266.489 & 254.632 & 0  & -  \\
  &   &  & 0.1 &  47.612& 47.641 & 434.338 & 0.190 & 0.210 & 137.920 & 68.101  & 0 & - \\
  &   &  & 1    & 47.612 & 49.375 & 365.626 & 0.558  & 0.190 & 325.569 & 63.173 & 0 & - \\
  \cmidrule(lr){3-13}
  &   & \multirow[t]{3}{*}{8} 
      & 0    & 62.326 & 62.326 & 528.232 & 0.053 & 0.370  & 1258.341 &  275.162 & 3 & 0.061 \\
  &   &  & 0.1 &  62.326 & 62.338 & 271.678 & 0.168 & 0.270 & 1904.136 & 2030.628 & 5 & 0.016  \\
  &   &  & 1    & 62.326 & 63.435  & 242.316 &   0.634 & 0.450  & 2361.074 & 3600  & 6 & 0.027  \\
\cmidrule(lr){2-13}
  & \multirow[t]{6}{*}{20}
      & \multirow[t]{3}{*}{4} 
       & 0    & 48.477 & 48.477 & 542.272 & 0.137 & 0.160 & 2585.092 & 3006.705 & 5 & 0.008 \\
  &   &  & 0.1 & 48.477 & 48.661 & 439.231 & 0.179 & 0.095 & 424.945 & 223.481 & 0 & - \\
  &   &  & 1    & 48.477 & 52.123 & 392.991 & 1.034 & 0.080 & 494.033 & 572.562 & 0  & - \\
  \cmidrule(lr){3-13}
  &   & \multirow[t]{3}{*}{8} 
      & 0    & 64.337 & 64.337  & 380.707 & 0.141 & 0.355 & 3600 & 3600 & 10 & 1.000 \\
  &   & & 0.1 & 64.337 & 64.400 & 173.389 & 0.166  & 0.135 & 3054.681 & 3600 & 8 & 0.037\\
  &   &  & 1    & 64.337 & 64.553 & 170.521  & 0.225 & 0.195 & 3011.619 & 3600  & 7 & 0.058 \\
\midrule
\multirow[t]{12}{*}{200}  
  & \multirow[t]{6}{*}{10}
      & \multirow[t]{3}{*}{4} 
       & 0    & 94.798 & 94.798 & 1494.580 & 0.129 & 0.300 & 1641.534 & 1217.67 & 1 &  0.000\\
  &   &  & 0.1 & 94.798 & 94.798 & 928.793 & 0.162 & 0.190 & 744.629 & 354.583 & 0 & - \\
  &   &  & 1    & 94.798 & 98.313 & 784.095 & 0.529 & 0.160 & 990.485 & 629.348 & 1 &  0.002\\
  \cmidrule(lr){3-13}
  &   & \multirow[t]{3}{*}{8} 
      & 0    & 124.591  & 124.591 & 843.438 & 0.053 & 0.290 & 2747.194  & 3600 & 6 & 0.125 \\
  &   & & 0.1 & 124.591 & 125.555 &  353.156 & 0.325 &  0.260 & 1744.012 & 1468.679 & 3 & 0.018  \\
  &   &  & 1    & 124.591 &  125.223 & 316.732 & 0.255 &  0.260 & 2900.122 & 3600 & 6 & 0.023 \\
\cmidrule(lr){2-13}
  & \multirow[t]{6}{*}{20}
      & \multirow[t]{3}{*}{4} 
       & 0    & 96.758 & 96.758  & 1080.307 & 0.134 & 0.185 & 3600 & 3600 & 10 & 0.804 \\
  &   &  & 0.1 & 96.758 & 97.042 & 843.456  & 0.166 & 0.085 & 943.262 & 444.74  & 1 & 0.007\\
  &   &  & 1    & 96.758 & 101.654  & 786.017 & 0.614 & 0.105 & 1786.051 & 1771.548  & 2 & 0.456\\
  \cmidrule(lr){3-13}
  &   & \multirow[t]{3}{*}{8} 
      & 0    & 128.686 & 128.701 & 833.982 & 0.173 & 0.075 & 3600 & 3600 & 10 &  1.000\\
  &   & & 0.1 &  128.686 & 128.778 & 326.169 & 0.129  & 0.105 & 3078.137 & 3600 & 8 & 0.121\\
  &   &  & 1    & 128.686 & 131.477 & 315.777 & 0.712  & 0.205  & 3600 & 3600 & 10 & 0.215\\
\bottomrule
\end{tabular}
\end{adjustbox}

\smallskip
{\footnotesize All reported values are averages computed over ten randomly generated instances. Maximum computational time is 3,600 seconds. The column TL shows the number of instances (out of 10) that reached the time limit. The column Gap reports the average relative MIP optimality gap for those instances that reached the time limit; otherwise, it is marked as `-`. The number of competitors was fixed to five ($|E|=5$) for all instances.}
\end{table}

\begin{table}[t]
\centering
\caption{Average construction time for the initial feasible solution and the Wasserstein bound calculation}
\label{tab:initial_times}
\begin{tabular}{ccc|cc}
\toprule
$N$ & $D$ & $B$ & Warm Start Time [s] & Bound Time [s] \\
\midrule
\multirow[t]{4}{*}{100} 
  & \multirow[t]{2}{*}{10} 
    & 4 & 0.105 & 0.505 \\
  &   & 8 & 0.107 & 0.852 \\
  \cmidrule(lr){2-5}
  & \multirow[t]{2}{*}{20} 
    & 4 & 0.249 & 1.271 \\
  &   & 8 & 0.243 & 1.414 \\
\midrule
\multirow[t]{2}{*}{200} 
  & \multirow[t]{2}{*}{10} 
    & 4 & 0.189 & 1.355 \\
  &   & 8 & 0.195 & 1.932 \\
  \cmidrule(lr){2-5}
  & \multirow[t]{2}{*}{20} 
    & 4 & 0.502 & 4.248 \\
  &   & 8 & 0.475 & 3.938 \\
\bottomrule
\end{tabular}

\smallskip
{\footnotesize All times are reported in seconds and represent averages across the test instances.}
\end{table}

We report the results in Table~\ref{tab:results}. They reveal several findings regarding both the computational performance of the formulation and the structural properties of the resulting explanations (captured demand, sparsity, and distributional shift). First, as expected, the captured demand is consistently maintained or slightly increased across all configurations. Thus, we confirm that the relative performance constraint is effectively enforced.

In terms of computational efficiency, the effect of the regularization parameter $\lambda$ is clearly visible. The unregularized model ($\lambda=0$) consistently exhibits longer average and median computing times compared to models with moderate regularization (e.g., $\lambda=0.1$). This effect becomes more pronounced as the instance size grows (larger values of $N$, $D$, and $B$). In addition, instances that reach the time limit tend to exhibit larger optimality gaps in the unregularized case.

Adding moderate Wasserstein regularization also improves the sparsity of the solutions: it reduces the complexity, i.e., the percentage of covariates changed (the $\ell_0$ ``norm" of the change). Moderate regularization encourages the model to concentrate adjustments in fewer facilities, rather than distributing small modifications across many, thereby leading to more interpretable explanations without substantially increasing the $\ell_1$ norm. Stronger regularization ($\lambda=1$) does not consistently bring these benefits and can even increase the percentage of covariates changed (see Table~\ref{tab:results}); we report this as an empirical observation.

In addition, the Wasserstein regularization promotes \emph{plausible}  explanations in a distributional sense: by penalizing large shifts between the factual and counterfactual demand distributions, it ensures that the induced choice behavior remains close to the observed (factual) demand pattern.

From a scalability perspective, the difficulty of the counterfactual problem largely mirrors that of the underlying factual CFLP. In particular, instances with larger budgets $B$ are significantly harder to solve and frequently reach the time limit, while increasing the number of customers $N$ and facilities $D$ also leads to higher runtimes. This aligns with the runtime behavior reported for the factual problem \citep[see][]{haasecomparison}.

Beyond a moderate level, regularization shows diminishing and eventually negative returns. Most of the reduction in $\mathcal{W}_2^2$ is already achieved at $\lambda=0.1$; increasing $\lambda$ to 1 yields only a small additional distributional gain (Table~\ref{tab:results}), while the $\ell_1$ cost generally rises and the problem becomes harder to solve, with more instances reaching the time limit and larger optimality gaps than at $\lambda=0.1$. In these results, strong regularization thus adds covariate cost and computational burden without a commensurate distributional benefit, which further supports a moderate choice. Across a finer grid of $\lambda$, the explanations remained stable around $\lambda=0.1$, which sits at the knee of this trade-off, removing most of the distributional shift at a small $\ell_1$ cost.

\subsection{Case Study: Montreal Electric Vehicle Charging Network}\label{sec:casestudy}

In this section, we compute relative counterfactual explanation for a real-world infrastructure planning problem in Montreal. We evaluate the model's ability to provide actionable insights for planners managing the Hydro-Qu\'ebec Circuit \'Electrique charging network.

\subsubsection{Data and Model Specification}

We build on the dataset and discrete choice models estimated by \citet{lamontagne2025makes}, who use preference data from charging sessions in Montreal to estimate MNL and mixed logit models for EV charging station selection. We adopt their MNL specification, consistent with the choice-based CFLP formulation of Section~\ref{sec:CFL}. We study two Montreal neighborhoods that present contrasting urban profiles: \textbf{Le Sud-Ouest \& Verdun}, a socio-economically diverse area that acts as a major 
transportation hub, including the Turcot Interchange and the Lachine Canal 
corridor ($|N|=444$, $|D|=20$); and \textbf{Plateau-Mont-Royal}, a residential neighborhood with high foot traffic but comparatively fewer existing charging stations ($|N|=308$, $|D|=11$). This contrast allows us to examine how the explanations differ across settings with different levels of existing infrastructure competition.

\paragraph{Facilities.}
For each neighborhood, we consider the subset of the 200 candidate locations from \citet{lamontagne2025makes} that fall within the neighborhood boundaries: $|D|=20$ candidates in Le Sud-Ouest \& Verdun and $|D|=11$ in Plateau-Mont-Royal. The set of existing competitors $E$ is represented as a single aggregate competitor, so that the sum over $E$ in Section~\ref{sec:explanations-cflp} reduces to this single term, with $b_n = \max_{e \in E} \exp(v_{e}^n)$ capturing the best available existing Level-2 station for each user $n$, consistent with the data structure 
provided by \citet{lamontagne2025makes}.

\paragraph{Users.}
Customer locations are represented by the centroids of dissemination areas within each study region, weighted equally ($q_n = 1/|N|$).There are $|N|=444$ users for Le Sud-Ouest \& Verdun and $|N|=308$ users for Plateau-Mont-Royal.

\paragraph{Utility specification.}
Following \citet{lamontagne2025makes}, the deterministic utility for user $n$ choosing station $d$ follows a near/far-from-home specification. Since 99.87\% of customer-station pairs have a Euclidean distance exceeding 400 meters (the near-home threshold), we apply the far-from-home specification uniformly across all pairs, introducing a negligible approximation:
\begin{equation}
    v_{d}^n = \underbrace{\beta_{\text{ctx}}\, x_d}_{\bar v_d(x_d;\,\beta_x)}
            + \underbrace{\beta_{\text{distFar}}\, \text{dist}_{nd}
            + \sum_{k \neq \text{ctx}} \beta_k\, s_{kd}}_{\hat v_d^n(\theta;\,\beta_\theta)}
    \label{eq:utility_ev}
\end{equation}
This is the decomposition $v_d^n = \bar v_d(x_d;\beta_x) + \hat v_d^n(\theta;\beta_\theta)$ of Section~\ref{sec:explanations-cflp}: the modifiable component $\bar v_d$ depends on the single covariate $x_d$ through $\beta_x = \beta_{\text{ctx}}$, while the fixed component $\hat v_d^n$ collects the attributes in $\theta$ that cannot be altered---the network distance $\text{dist}_{nd}$ (km) between dissemination area $n$ and station $d$, and the remaining station attributes $s_{kd}$ (gas station indicator, log-density of fast food, supermarkets, shopping, malls, leisure and sport), with parameters $\beta_\theta = (\beta_{\text{distFar}}, \{\beta_k\}_{k\neq\text{ctx}})$.
Parameter values are taken from \citet{lamontagne2025makes} (L2-MNL model): $\beta_{\text{distFar}} = -0.054$, $\beta_{\text{outletsFar}} = 0.029$ (used as $\beta_{\text{ctx}}$ in Experiment~1), and $\beta_{\text{rest}} = 0.117$ (used as $\beta_{\text{ctx}}$ in Experiment~2).

\paragraph{Covariates.}
Following the framework of Section~\ref{sec:CFL}, we factorise $\exp(v_d^n) = \phi_d(x_d) \cdot a_{nd}$, where 
$\phi_d(x_d) = \exp(\beta_{\text{ctx}} \cdot x_d)$ captures the covariate component that can be changed and $a_{nd}$ is fixed throughout optimization. We conduct two experiments, each with a different choice of covariate $x_d$, which play conceptually distinct roles.

In \textbf{Experiment~1}, $x_d = \text{outlets}_d$ with 
$\beta_{\text{ctx}} = 0.029$, starting from $x_d^0 = 1$ outlet for all candidates. The number of outlets is an operational decision of the network operator rather than an exogenous feature. The counterfactual explanation thus provides \emph{actionable} insights: it identifies the minimal change in outlet capacity across candidate stations such that a solution including the target station performs at least as well as the factual solution in terms of captured demand. The result is interpreted as $\lceil x_d^{\text{new}} \rceil$ outlets in practice. Notably, the explanation may require adjusting outlets at stations other than the target station itself, reflecting the substitution patterns induced by the MNL model.

In \textbf{Experiment~2}, $x_d = \log(1 + n_{\text{rest},d})$ with $\beta_{\text{ctx}} = 0.117$, where $n_{\text{rest},d}$ is the number of restaurants within 300 metres of station $d$, with factual values taken from \citet{lamontagne2025makes}. Unlike Experiment~1, this covariate is an exogenous feature of the urban environment beyond the operator's direct control. The explanation identifies which changes in surrounding land use would render the target station worth opening, connecting the optimization to urban planning decisions such as zoning and commercial development. Results are back-transformed as $n_{\text{rest}}^{\text{new}} = 
\exp(x_d^{\text{new}}) - 1$.

Together, the two experiments demonstrate the flexibility of the framework: it can generate explanations in terms of both operator-controlled decisions and exogenous environmental factors.

\subsubsection{Counterfactual Query and Properties}

For each neighborhood, the factual solution $\bm z^0$ is obtained by solving problem~\eqref{eq:simplifyCFL} with budget $B=5$, using the full utility specification~\eqref{eq:utility_ev} evaluated at the factual covariate values from \citet{lamontagne2025makes}; it is therefore common to both experiments. The desired space $\mathcal{D}$ fixes a single target station $d^{*}$ to be open in the counterfactual solution. We select $d^{*}$ as the candidate closest to the neighborhood centroid that is absent from $\bm{z}^0$ and has no nearby restaurants. This choice is driven by a practical question a stakeholder might ask: a centrally located site seems like a natural candidate for planning. If it is excluded from the optimal solution, that decision calls for a clear explanation. In Le Sud-Ouest \& Verdun, $d^* = \texttt{110}$; in Plateau-Mont-Royal, $d^* = \texttt{68}$. 

We solve Problem~\eqref{eq:RE2} with $\alpha = 1.0$ and $\lambda = 0.1$. Rather than tune $\lambda$ separately for the case study, we fix $\lambda = 0.1$ a priori from the synthetic analysis of Section~\ref{sec:larger}. This transfer is meaningful because the Wasserstein term is normalized by the model-free reference value of Section~\ref{sec:lowerandwarmstart}, computed per instance, $\lambda$ is a relative weight between feature-level dissimilarity and distributional regularization with an instance-independent interpretation, so a value calibrated on one set of instances carries over to others.

The choice of $\epsilon=10^{-5}$ in~\eqref{eq:LB} satisfies the condition required for the lower bound (Section~\ref{sec:lowerandwarmstart}): it is well below the smallest positive choice probability observed across the factual and counterfactual solutions of all instances (almost four orders of magnitude above $\epsilon$).

The explanations are evaluated with respect to the desirable properties introduced before. \emph{Proximity} and \emph{sparsity} are directly induced by the $\ell_1$ term, which penalizes large and non sparse changes in $\phi_d$. \emph{Actionability} is naturally captured by Experiment~1: the  suggested change in number of outlets is a concrete, implementable  decision for the operator.  Regarding \emph{plausibility}, in the counterfactual explanation literature this typically requires that the suggested changes remain 
consistent with the data manifold, i.e., that the counterfactual instance could plausibly be observed in practice. In our setting, since only a single covariate is modified, the natural analogue is to bound the changes within the range of values observed across all candidate locations in the dataset. We therefore impose 
$x_d^{\text{new}} \leq x_{\max}$, where $x_{\max}$ corresponds to the maximum observed value in the Montreal dataset: at most 13 outlets and at most 39 nearby restaurants. This ensures that the suggested changes are plausible. Moreover, in a distributional sense, plausibility is additionally addressed through the Wasserstein regularization: by penalizing large shifts between the factual and counterfactual demand  distributions, we ensure that the induced choice behavior remains close to realistic patterns under the MNL model.

\subsubsection{Results}

\begin{figure}[!htbp] 
\centering
\begin{tabular}{cc}
    \includegraphics[width=0.45\textwidth]{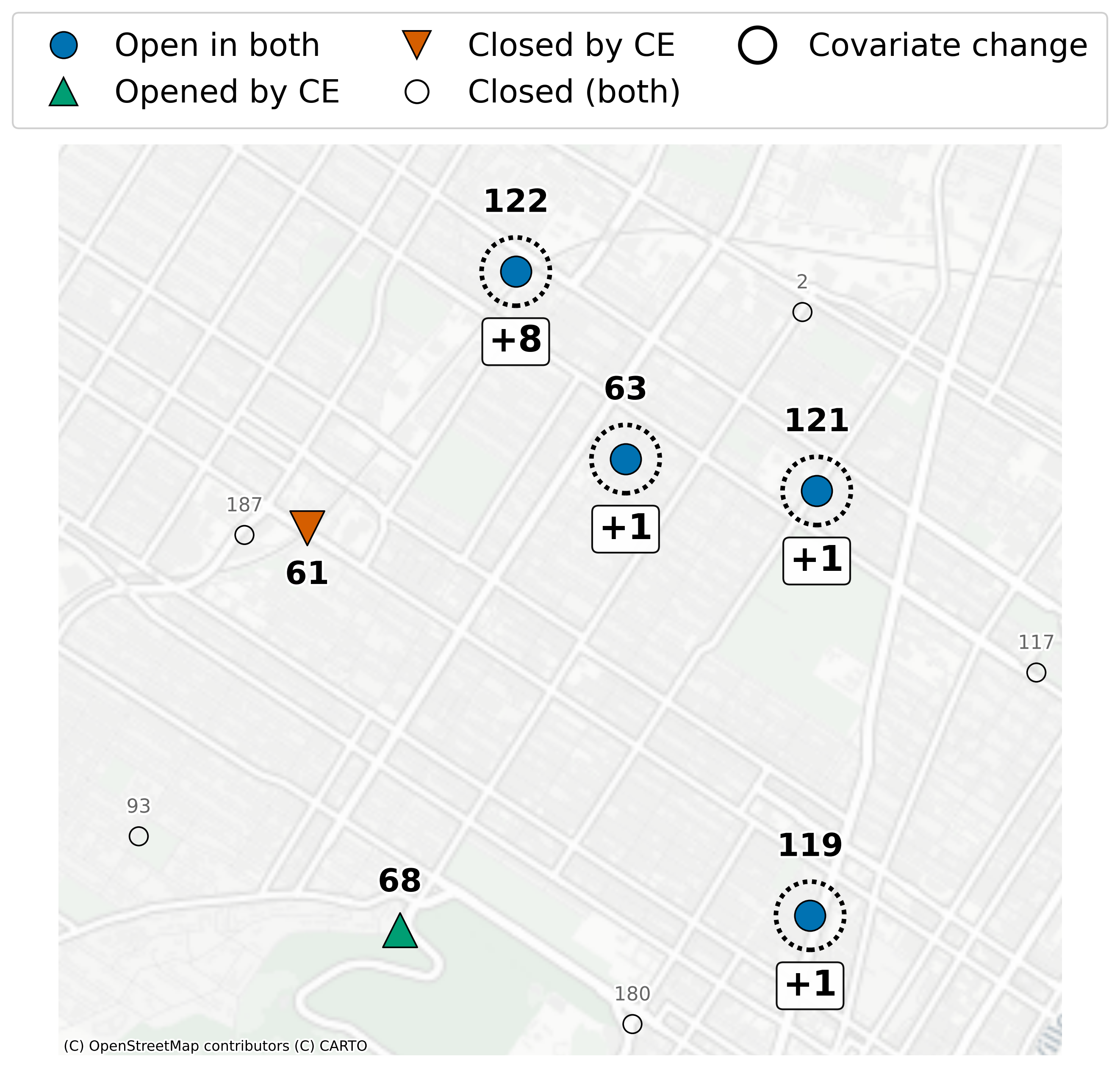} & 
    \includegraphics[width=0.45\textwidth]{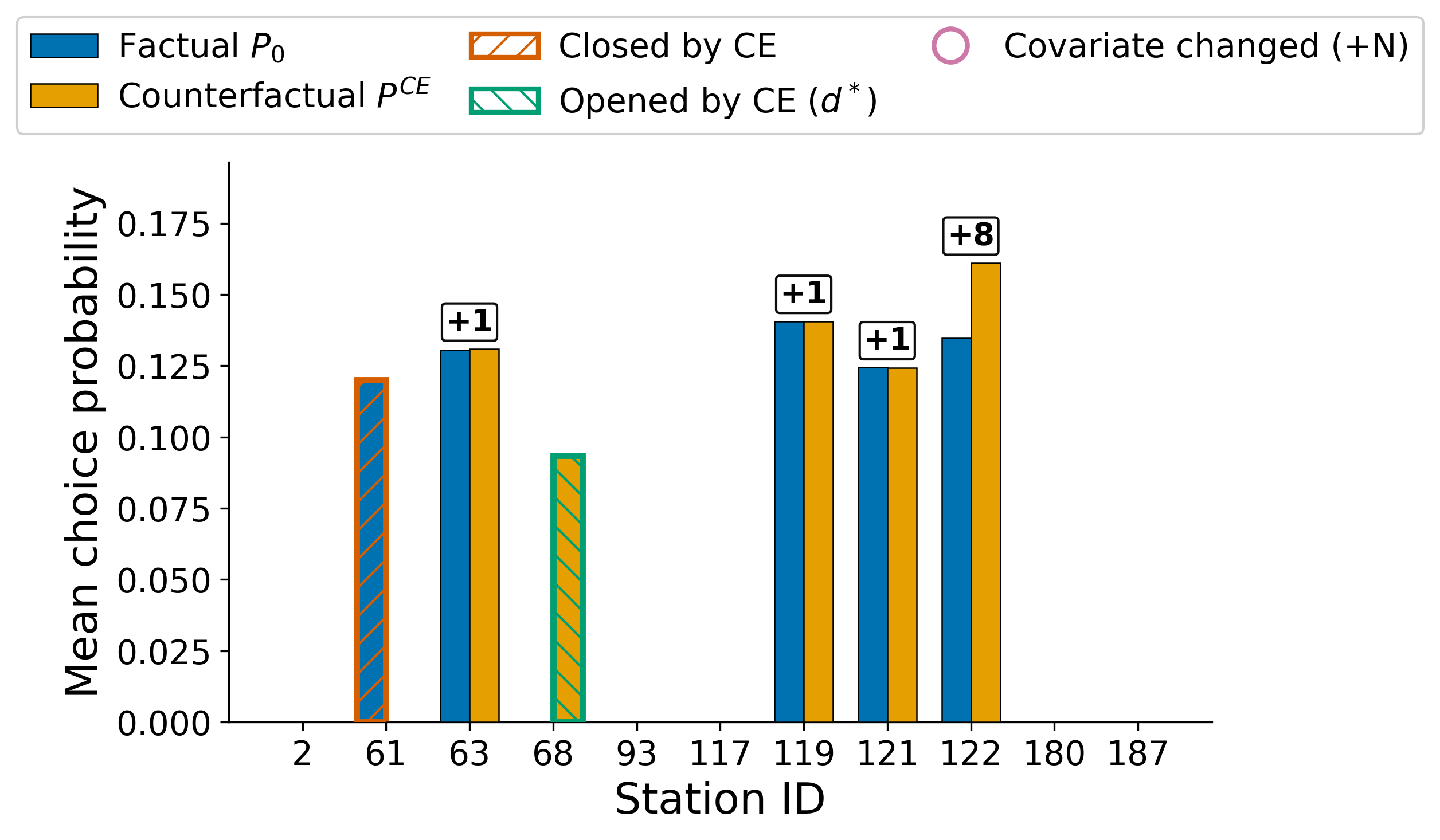} \\
    (a) Map: outlets & (b) Distributions: outlets \\ \addlinespace[10pt]
    \includegraphics[width=0.45\textwidth]{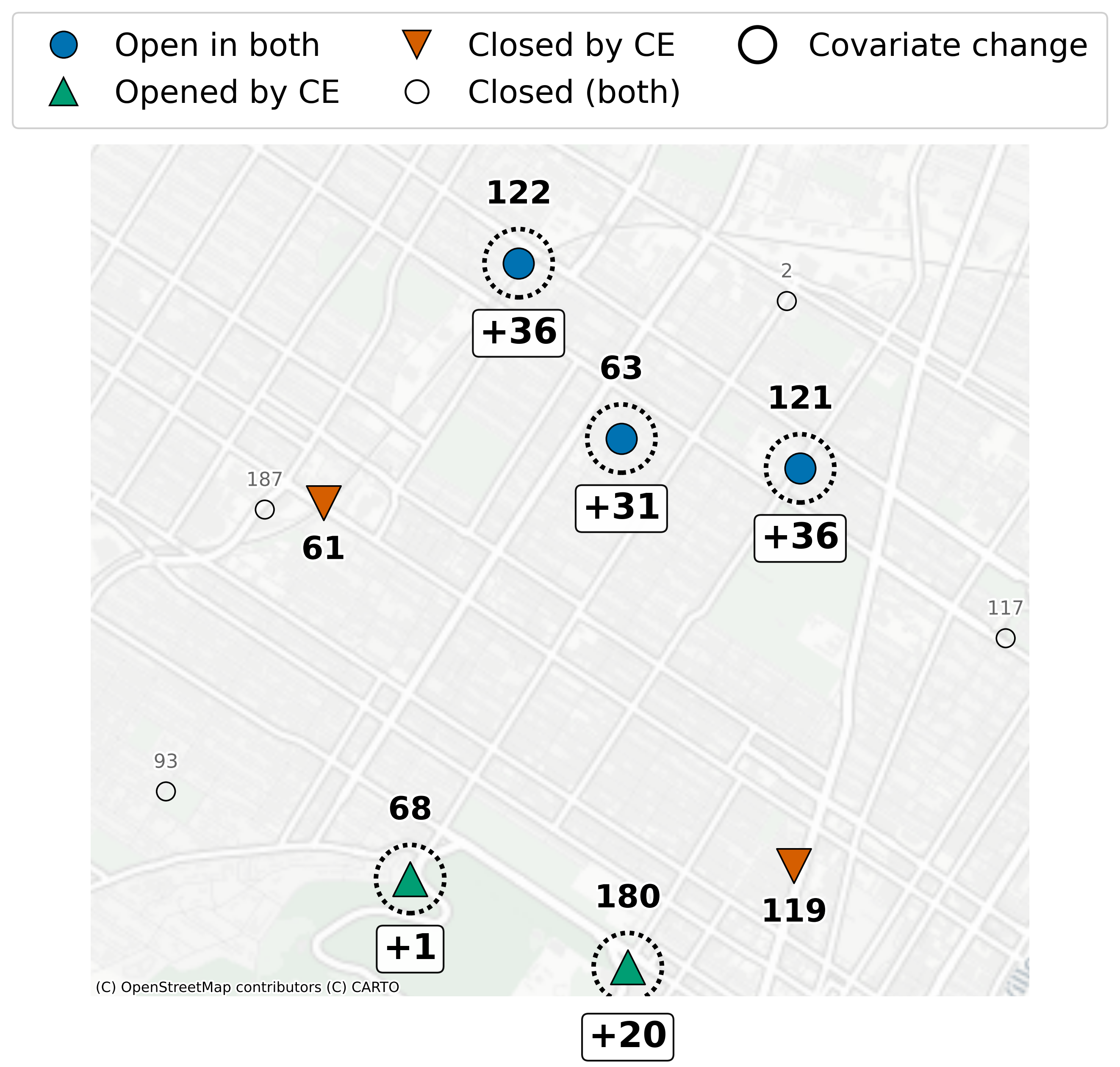} & 
    \includegraphics[width=0.45\textwidth]{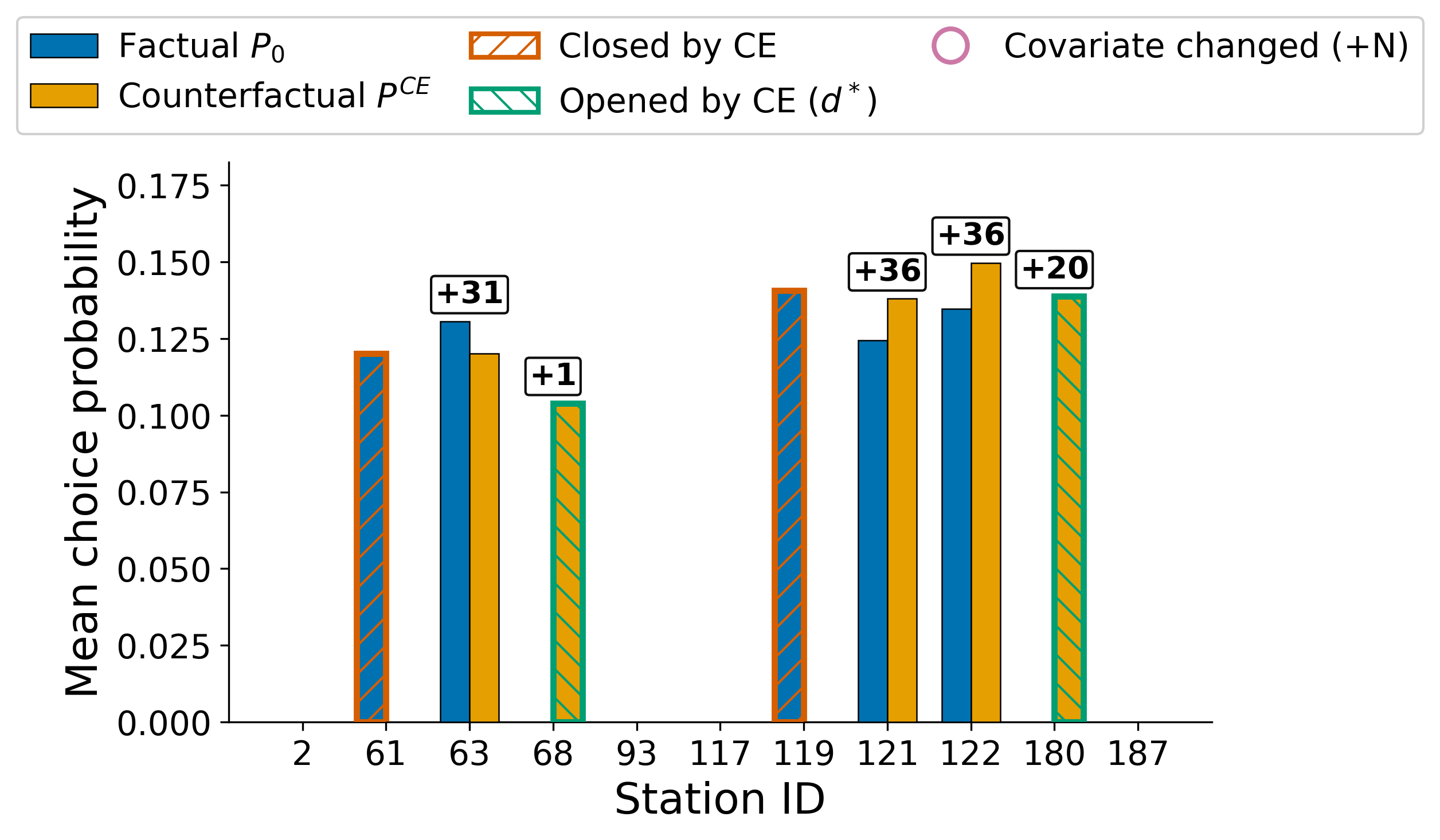} \\
    (c) Map: restaurant & (d) Distributions: restaurant
\end{tabular}
\caption{Counterfactual explanations for Plateau-Mont-Royal. Left: map of factual and counterfactual solutions. Right: mean choice probability per station under factual ($P_0$, blue) and counterfactual ($P^{CE}$, orange) solutions. 
In the maps, a blue circle denotes a station open in both solutions, a green upward triangle ($\triangle$) a station opened by the counterfactual, a vermilion downward triangle ($\triangledown$) a station closed by the counterfactual, and a small hollow circle a station closed in both. A dotted black ring around a station marks a covariate change, with the adjacent integer label giving the number of outlets (Experiment~1) or restaurants (Experiment~2) to be added. Each station is annotated with its identifier.}
\label{fig:plateau}
\end{figure}

\begin{figure}[!htbp]
\centering
\begin{tabular}{cc}
    \includegraphics[width=0.45\textwidth]{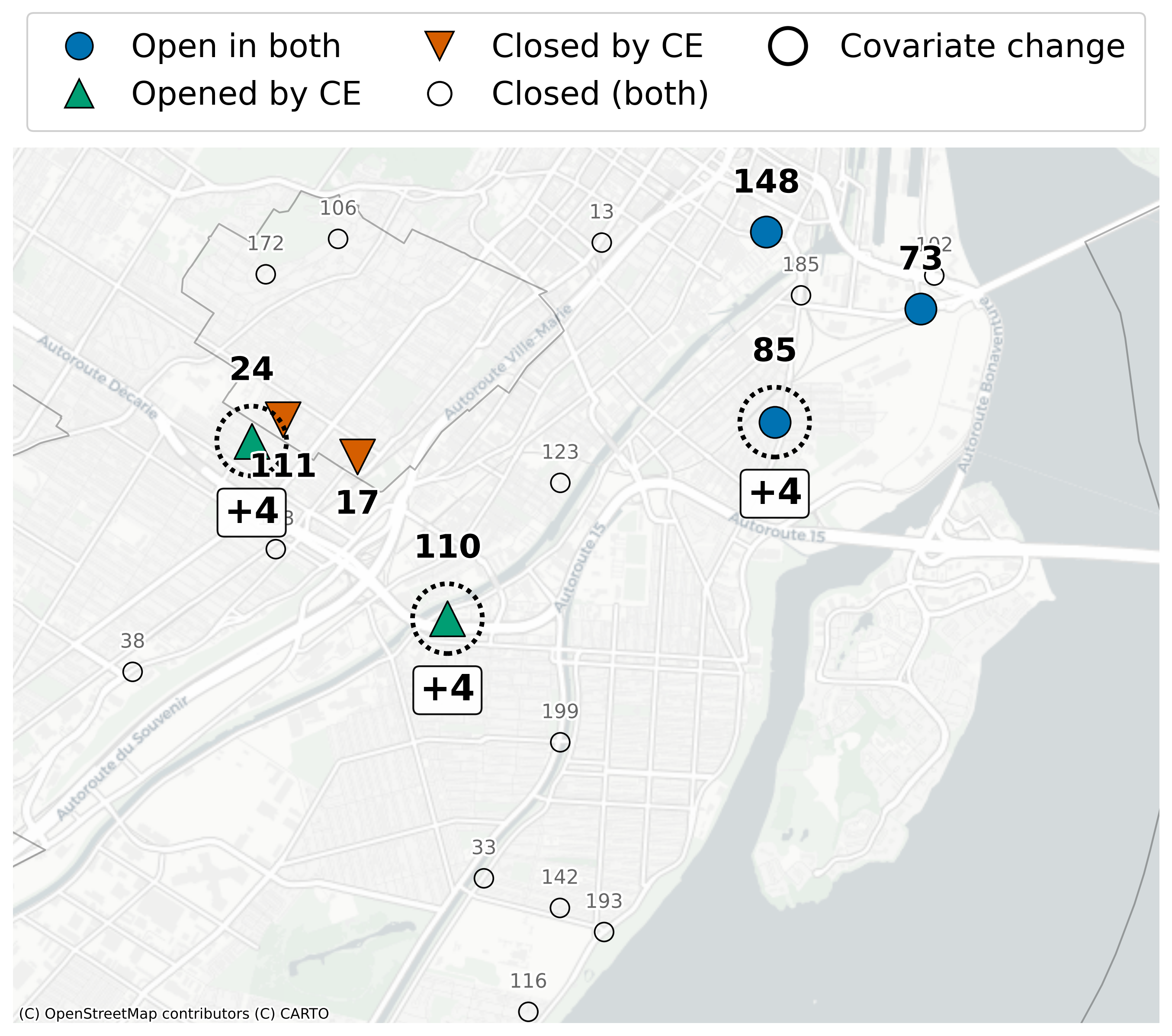} & 
    \includegraphics[width=0.45\textwidth]{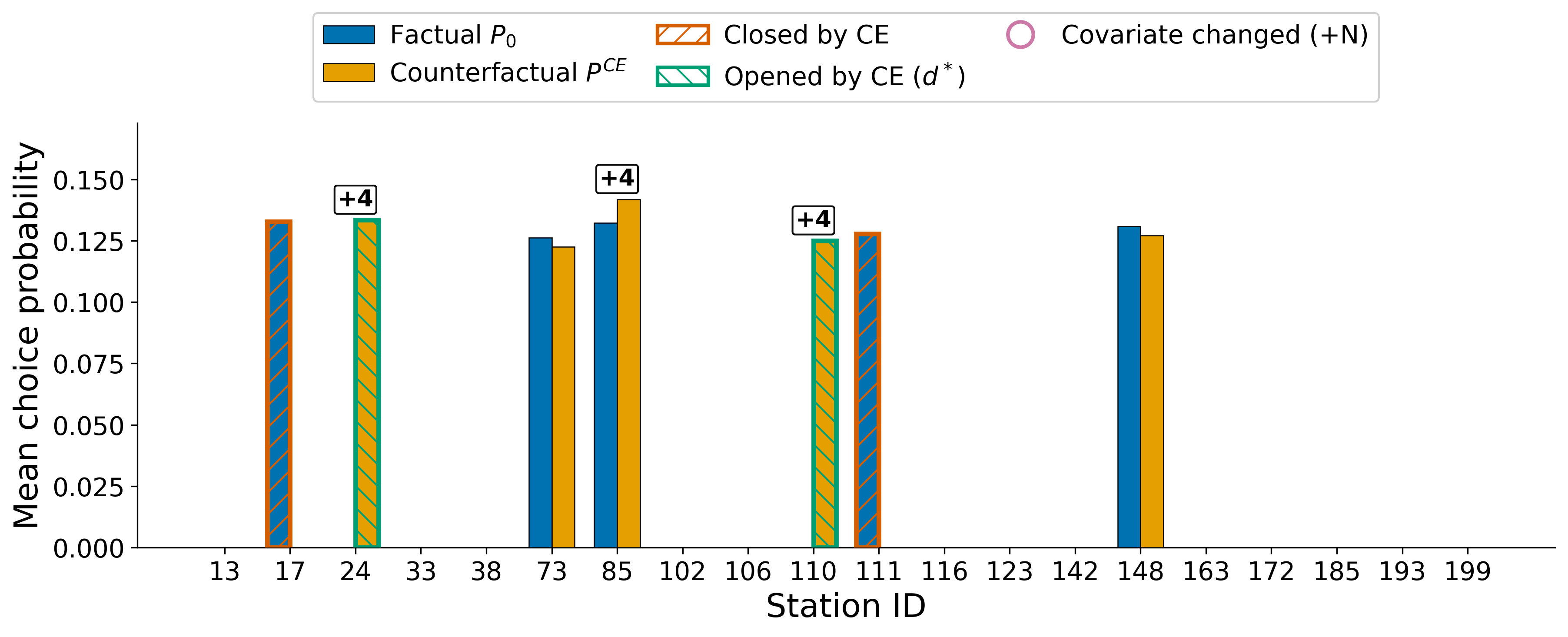} \\
    (a) Map: outlets & (b) Distributions: outlets \\ \addlinespace[10pt]
    \includegraphics[width=0.45\textwidth]{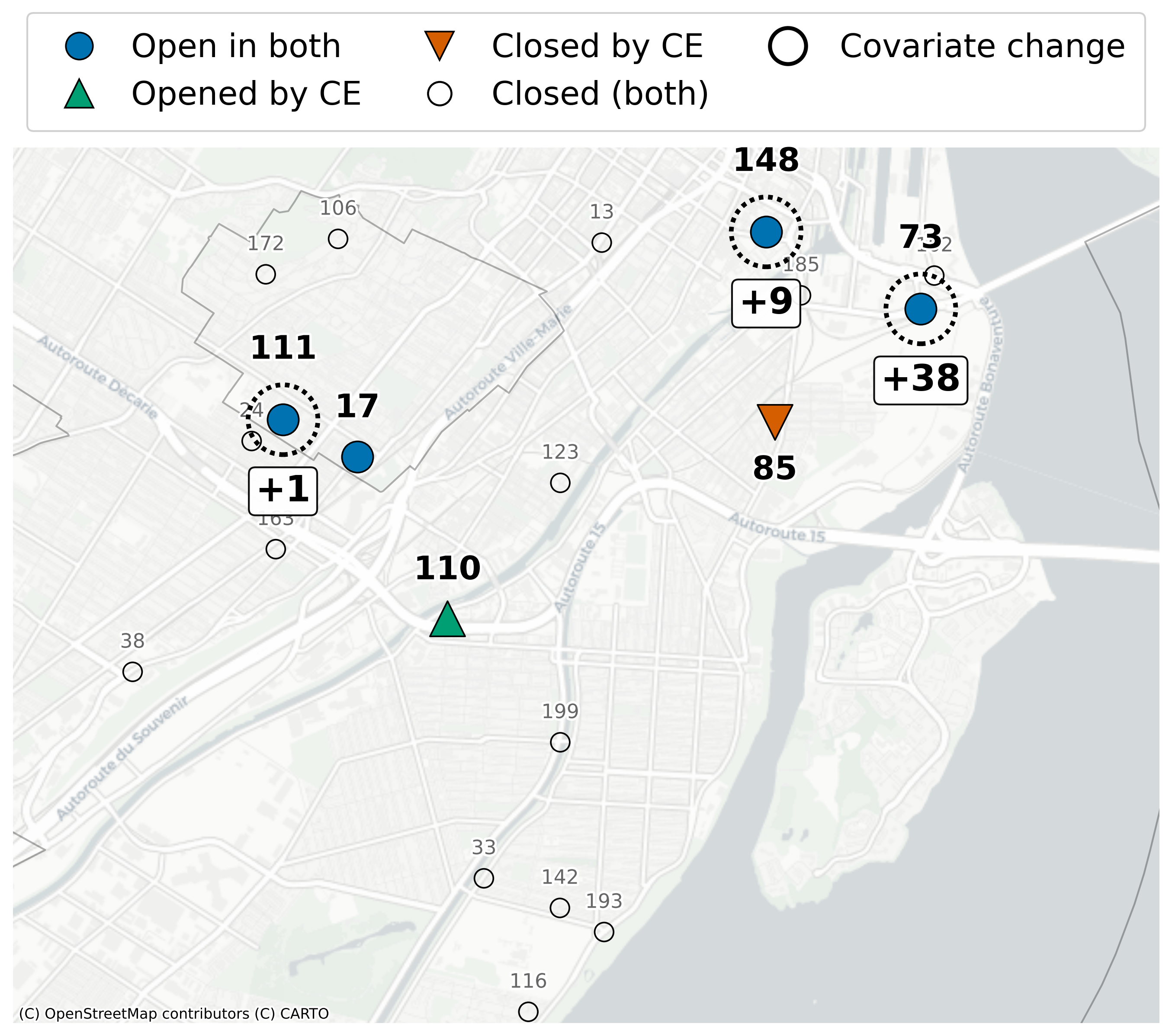} & 
    \includegraphics[width=0.45\textwidth]{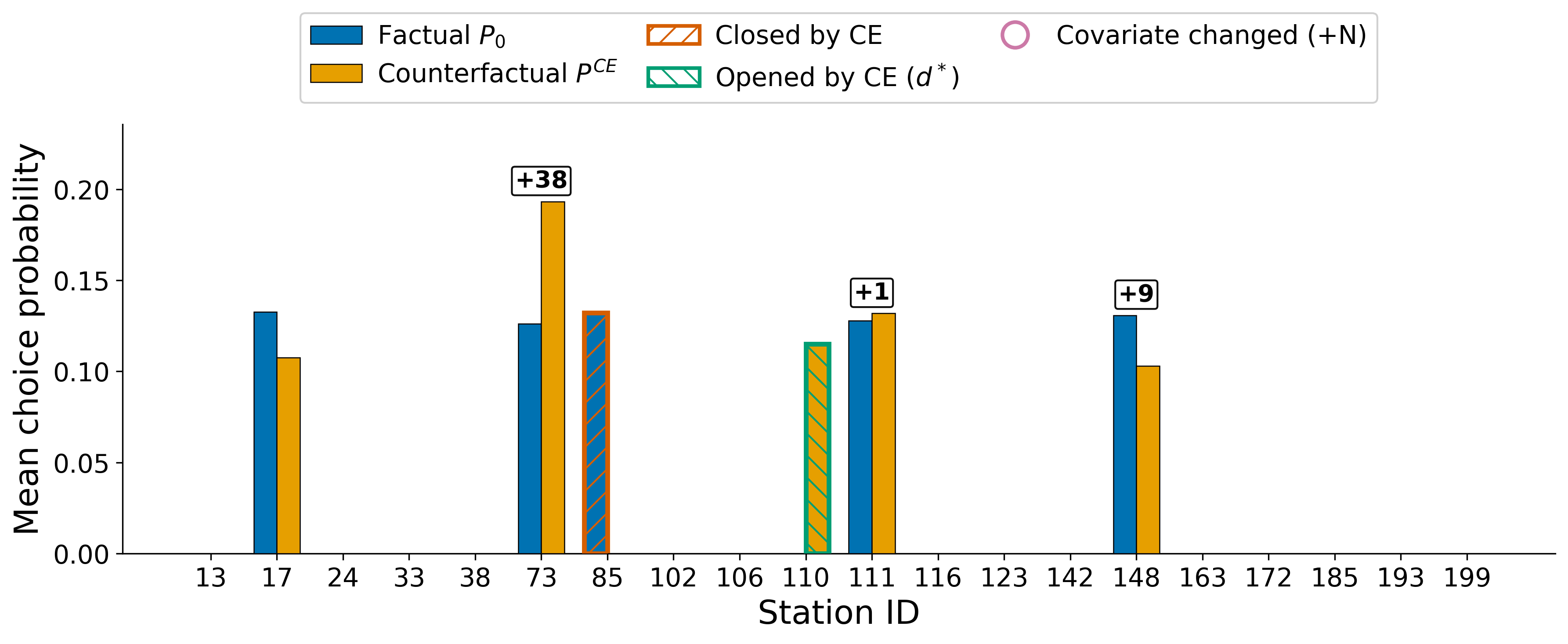} \\
    (c) Map: restaurant & (d) Distributions: restaurant
\end{tabular}
\caption{Counterfactual explanations for Le Sud-Ouest \& Verdun. Left: map of factual and counterfactual solutions. Right: mean choice probability per station under factual ($P_0$, blue) and counterfactual ($P^{CE}$, orange) solutions. In the maps, a blue circle denotes a station open in both solutions, a green upward triangle ($\triangle$) a station opened by the counterfactual, a vermilion downward triangle ($\triangledown$) a station closed by the counterfactual, and a small hollow circle a station closed in both. A dotted black ring around a station marks a covariate change, with the adjacent integer label giving the number of outlets (Experiment~1) or restaurants (Experiment~2) to be added. Each station is annotated with its identifier.}
\label{fig:sudouest}
\end{figure}

We report the results for both neighborhoods and experiments in Figures~\ref{fig:plateau} and~\ref{fig:sudouest}. For each experiment, we show the map of the factual and counterfactual solutions and the induced choice probability distributions aggregated over customers.

\paragraph{Plateau-Mont-Royal.}
The factual solution opens stations \texttt{61}, \texttt{63}, \texttt{119}, \texttt{121}, and \texttt{122}, capturing $Q^{\text{factual}} = 0.6504$, reported to four decimals to distinguish it from Le Sud-Ouest \& Verdun ($0.6498$), which agrees to three. The desired station $d^* = \texttt{68}$ is excluded.

In Experiment~1, the counterfactual closes \texttt{61} and opens \texttt{68}, with outlet changes at four stations: \texttt{122} requires 9 outlets, \texttt{63}, \texttt{119} and \texttt{121} require 2 outlets each. As shown in Figure~\ref{fig:plateau}(b), the demand distribution shifts as follows: the probability mass of station \texttt{61}  ($P_0 = 0.120$), which closes, is redistributed primarily to the newly opened \texttt{68} ($P^{CE} = 0.094$) and partially to  \texttt{122} ($P_0 = 0.135 \to P^{CE} = 0.161$), whose increased  outlet capacity absorbs the remaining displaced demand. The remaining open stations are unaffected. The Wasserstein regularization promotes this distributional stability, which gives plausible explanations.

In Experiment~2, the counterfactual closes \texttt{61} and \texttt{119}, and opens \texttt{68} and \texttt{180}. The explanation requires three stations to reach their maximum observed restaurant density (39 restaurants), with \texttt{68} itself requiring at least 1 nearby restaurant. The required environmental changes are substantially larger than in Experiment~1, suggesting that the commercial environment of Plateau-Mont-Royal does not naturally support the inclusion of \texttt{68} without significant surrounding amenity improvement.

\paragraph{Le Sud-Ouest \& Verdun.}
The factual solution opens stations \texttt{17}, \texttt{73}, \texttt{85}, \texttt{111}, and \texttt{148}, capturing $Q^{\text{factual}} = 0.6498$. The desired station $d^* = \texttt{110}$ is excluded.

In Experiment~1, the counterfactual closes \texttt{17} and \texttt{111}, and opens \texttt{110} and \texttt{24}. Three stations require outlet increases to 5. The demand distribution (Figure~\ref{fig:sudouest}(b)) shifts 
as follows: the probability mass of \texttt{17} ($P_0 = 0.133$) and \texttt{111} ($P_0 = 0.128$), which close, is absorbed by the newly opened \texttt{110} ($P^{CE} = 0.125$) and \texttt{24} ($P^{CE} = 0.133$), with \texttt{85} gaining slightly ($P_0 = 0.126 \to P^{CE} = 0.141$).

In Experiment~2, the counterfactual closes \texttt{85} and opens \texttt{110}. Station \texttt{73} requires 38 nearby restaurants (currently 0) and \texttt{148} increases from 12 to 21. As visible in  Figure~\ref{fig:sudouest}(d), the shift in demand towards \texttt{73} is pronounced as its mean probability increases from $0.13$ to $0.19$, reflecting the large covariate change required at that station.

\paragraph{Discussion.}
Across both neighborhoods, we can see a clear asymmetry. Outlet-based explanations are  sparse and the induced distributional change is localized: the probability mass of the closing station shifts mainly to the newly opened $d^{*}$ (and to a nearby station absorbing the remainder), while the other open stations remain largely unaffected. This reflects that moderate capacity investments suffice to justify including $d^{*}$. Restaurant-based explanations require substantially larger  environmental changes, frequently reaching the maximum observed  values, suggesting that the exclusion of these stations is driven  primarily by their unfavorable commercial environment rather than
insufficient infrastructure capacity.

More broadly, the case study illustrates a complementary role for counterfactual explanations alongside traditional sensitivity analysis: while sensitivity analysis measures how the optimum responds to perturbations of the inputs, a relative counterfactual answers the inverse question of how much the inputs would need to change for a specific decision to be supported by the model. This makes the framework an ex-ante diagnostic, valuable when a planner holds a strong prior about a particular station before any infrastructure commitment. The diagnostic can take two qualitatively different forms, both visible in our experiments. In Experiment~1, it materializes as a concrete and proportionate intervention---a moderate redistribution of charging outlets across neighboring stations---signaling that the prior is consistent with the model up to a tractable adjustment that can guide actionable capacity planning. In Experiment~2, by contrast, the required changes are large and saturate the observed range, with up to 39 nearby restaurants at multiple stations, indicating that the gap between the planner's prior and the model's view is too large to bridge under realistic urban conditions. In both cases, the counterfactual converts a stakeholder prior into a quantitative target, giving them a concrete input either to act on or to revise.

\section{Conclusion}\label{sec:conclusion}

This paper introduces an approach for computing counterfactual explanations for contextual stochastic problems with endogenous uncertainty and binary decision variables. We propose the use of a Wasserstein regularization term and a model-free variant to compute a lower bound.

We apply our methodology to a choice-based CFLP under the MNL model, demonstrating that obtaining counterfactual explanations can be formulated as a mixed integer bilinear program. Including a moderate Wasserstein regularization term in the objective function reduces the computing time compared to an objective without regularization. Moreover, the resulting counterfactuals are sparser and lead to more gradual changes in the choice probabilities.

We further demonstrate the practical applicability of the framework through a real-world case study on electric vehicle charging station planning in Montreal, where the explanations provide actionable capacity recommendations for network operators and environmental insights for urban planners.

We conclude this paper by discussing the limitations of this work that warrant future research. First, like \citet{viviercf}, the current formulation requires the contextual covariates to be continuous. Extending the methodology to handle categorical covariates remains an open challenge, as additional constraints and a reformulation would be needed to explicitly encode their categorical structure.

Second, the solution is generally not unique, which can lead to potential disagreements among stakeholders and the possibility of manipulating explanations \citep{brughmans2024disagreement}. The existence of multiple valid counterfactual explanations may also raise ethical concerns, for example, by enabling the concealing of sensitive covariates. A crucial direction for future research is the incorporation of objective terms that ensure solution uniqueness.

Third, our application to the choice-based CFLP relies on the MNL model, whose independence of irrelevant alternatives (IIA) property is a well-known restriction. \citet{lamontagne2025makes} also report results for mixed logit models. Extensions to more flexible RUM models, such as mixed logit or generalized extreme value \citep{dam2022submodularity}, would relax IIA at the cost of additional algorithmic complexity, and are a natural direction for future work.

Fourth, our experiments indicate that the benefit of distributional regularization is non-monotone in $\lambda$: moderate values improve sparsity, plausibility, and tractability, whereas strong regularization erodes these gains, inflating both the covariate cost and the computational effort without a commensurate reduction in distributional shift. A principled, possibly instance-adaptive criterion for selecting $\lambda$, together with a fuller characterization of this over-regularization regime, is left for future work.

Other promising research avenues include the computation of weak instead of relative counterfactuals, where the optimality of the new solution is enforced, extensions to other discrete choice models, as well as applications in contextual problems such as assortment optimization and pricing.

\bibliographystyle{apalike}
\bibliography{bibliography}

\appendix

\section{Formulation from \cite{haase2009discrete}} \label{sec:HaaseForm}

\setcounter{equation}{0}
\renewcommand{\theequation}{A.\arabic{equation}}

Let $w_{nd}$ and $\hat w_n$ be non-negative decision variables. Then, \citet{haase2009discrete} show that \eqref{eq:simplifyCFL} can be reformulated as follows:

\begin{align}
    \max \quad&\sum_n q_n \sum_d w_{nd} \label{eq:h1}\\
    \text{s.t.} \quad & \hat w_n + \sum_d{w_{nd}} \leq 1 & n\in N&\label{eq:h2}\\
    & a_{nd}(w_{nd}-z_d) + b_n w_{nd} \leq 0 & d \in D \quad n \in N &\label{eq:h3}\\
    &w_{nd} - \frac{a_{nd}}{b_n} \hat{w}_n \leq 0 & d\in D \quad n \in N&\label{eq:h4}\\
    &\sum_d z_d =B & d\in D&\label{eq:h5}\\
    &z_d \in \{0,1\} & d\in D&\label{eq:h6}\\
    &w_{nd}\geq0 & d\in D \quad n \in N&\label{eq:h7}\\
    &\hat w_{n} \geq 0 & n\in N\label{eq:h8}
\end{align}

It is easy to see that $w_{nd}=\Prob_n(d)$ and $\hat w_n=\sum_{e\in E} \Prob_n(e)$.  Constraint~\eqref{eq:h2} ensures that the total probability assigned to each customer does not exceed one. Constraint~\eqref{eq:h4} reproduces the proportions given by the logit model, while constraint \eqref{eq:h3} is just a tighter bound than simply using $w_{nd}\leq z_d$. Constraint ~\eqref{eq:h5} defines the budget. Finally, constraints ~\eqref{eq:h6}-\eqref{eq:h8} specify the nature of the decision variables.

\end{document}